\documentclass[11pt,reqno]{amsart}
\usepackage{ifthen, amsmath,amssymb,graphicx,epstopdf,mathrsfs, amsfonts}
\usepackage{esint,bm}
\usepackage[italian,english]{babel}
\usepackage{booktabs}
\usepackage{subfig}
\usepackage{url}

\newenvironment{sistema}%
{\left\lbrace\begin{array}{@{}l@{}}}%
{\end{array}\right.}

\provideboolean{shownotes}
\setboolean{shownotes}{true} 
\newcommand{\margnote}[1]{
	\ifthenelse{\boolean{shownotes}}
	{\marginpar{\raggedright\tiny\texttt{#1}}}
	{}
	}

\newcounter{cont}[section]

\numberwithin{equation}{section}

\begin{document}
\baselineskip=16pt

\title[Burridge-Knopoff]{Burridge-Knopoff}

\author[P. Moschetta]{Pierfrancesco Moschetta}
\address[Pierfrancesco Moschetta]{Dipartimento di Matematica, Sapienza Universit\`a di Roma (Italy)}
	\email{pierfrancesco.moschetta@gmail.com}
 
\author[C. Mascia]{Corrado Mascia}
\address[Corrado Mascia]{Dipartimento di Matematica,
	Sapienza Universit\`a di Roma (Italy)}
	\email{mascia@mat.uniroma1.it}
 
\thanks{}


\begin{center}\sf
{\Large Assessment of Predictor-Corrector strategy 

for the Burridge-Knopoff model}\vskip.25cm

Pierfrancesco MOSCHETTA,  Corrado MASCIA
\end{center}
\vskip.5cm

\begin{quote}\footnotesize\baselineskip 12pt 
{\sf {\bf Abstract.}
A Predictor-Corrector strategy is employed for the numerical simulation of the one-dimensional 
Burridge-Knopoff model of earthquakes.
This approach is totally explicit and allows to reproduce the main features of the model.
The results achieved are compared with those of several previous works available in the literature,
in order to state the effectiveness of the novel numerical strategy.
Simulations are performed starting from the simplest cases and are aimed at studying the qualitative trends
of the phenomena under analysis.
By increasing the size of the associated differential system, it is possible to examine data on the basis of
the Gutenberg-Richter statistical law.
Finally, some tests are conducted to investigate the continuum limit of the discrete Burridge-Knopoff model
towards a macroscopic dynamics.} 
\end{quote}

\section{Introduction}\label{sec:intro}
\selectlanguage{english}
Earthquakes are doubtless an open research field. The necessity of improving our knowledge of this kind of geophysical mechanisms and related topics is very strong. 

Specifically, an earthquake occurs along fractures in the Earth's crust, named \textit{faults}, characterized by a steady accumulation of tension, when big quantities of energy are suddenly released due to the relative motion of the two sides involved. To understand better we recall that Earth's lithosphere includes the crust and is also composed of a part of the upper mantle; Moreover, it presents a complex structure divided into distinct blocks, the \textit{tectonic plates}. This point is crucial, because it is along the borders of tectonic plates that the great accumulation of tension we mentioned above takes place. The plates are indeed continually stressed by external forces, whose nature has been investigated and is continuing being object of study: scientists think this stress to be caused by the mantle convection but a gravitational contribution is not ignored \cite{Doglioni15}. 

A central role in our analysis is played by the friction. Indeed, although the existence of forces able to solicit plates is of course an important factor to explain seismic events, nothing would happen if friction did not inhibit the relative motion between the two different sides of an active fault. Strongly connected with these concepts is the \textit{stick-slip phenomenon}, firstly associated to the earthquakes by Brace and Byerlee \cite{Brace1966}. The borders of a fault exhibit asperities which make the local slip very difficult: as a consequence tension increases and the motion is inhibited by the balance between tension and friction. Once that this equilibrium is compromised, due to the steady accumulation of stress, a slip of the sides involved occur and a great quantity of energy is released, generating an earthquake. The alternation between period of quiescence, in which tension increases, and phases in which tension is released along the fault, through the motion of the plates, is a typical example of the stick-slip behaviour. 

It is important to notice that in the last decades a great effort has been made to investigate the statistical properties of earthquakes. This way of thinking is strongly connected with the idea of \textit{self-organized criticality}, SOC, developed by Bak et al. \cite{Bak1987} and its influence on seismic events \cite{Bak1989}. By following this concept lots of natural phenomena are explainable in terms of criticality: these kinds of processes can self-organize and reach critical states. When similar states are reached, little perturbations affecting the elements belonging to the systems can propagate and involve items of any size \cite{Winslow1997}. In the SOC view this behaviour is often illustrated by basic laws collecting the statistical properties of the process studied. As concerns the earthquakes, two important power laws would be a concrete manifestation of the SOC principles: the Gutenberg-Richter law \cite{GutRich56} for the magnitude distributions and the Omori law for the aftershocks sequences. In order to point out the SOC idea, the earthquake models are often analyzed with cellular automation approach, as in the work by Olami, Feder and Christensen \cite{Olami1992}.

One of the most famous mathematical models developed to study earthquakes and its statistical properties, especially pursuing the idea of a qualitative comparison with real phenomena, is the Burridge-Knopoff model, proposed by Robert Burridge and Leon Knopoff in 1967 \cite{BurrKnop67}. This is a deterministic dynamical system whose computational investigation provides lots of useful results to achieve a sufficiently accurate analysis of seismic events. The Burridge-Knopoff model has been deeply investigated in order to pursue a statistical study of earthquakes \cite{Kawamura2012} and continues to be a landmark in this research field, due to its nontriviality but, at the same time, its semplicity \cite{deArcangelis2016}. On a mathematical level the associated differential system exhibits a discontinuous right hand side, arising from the choice of the friction law. This is a direct consequence of the alternation caused by the stick-slip dynamics and expressed by the friction, the only source of nonlinearity in the model.
This alternation produces a dry friction. Lots of models arising from applications exhibit similar characteristics and require careful analysis. 
An analytical study of a non-smooth friction-oscillator model, qualitatively very close to the Burridge-Knopoff model, is provided in \cite{KunzKupp97}. Obviously also the related numerical problem must be adequately approached: in this sense some numerical methods are employed for non-smooth systems \cite{Acary2008} and suitable regularizations are often performed \cite{Guglielmi2015}.

In this paper we describe our numerical adjustment of the system and discuss the main results provided by simulations. Our purpose consists in employing a numerical method based on a Predictor-Corrector approach. 
The idea behind a Predictor-Corrector strategy \cite{Quarteroni} is inspired by the necessity of furnishing a good initial guess to start fixed-point iterations when an implicit method is invoked. Indeed, because several function evaluations are generally needed by using the fixed-point method, trying to reduce the computational cost becomes important. So the basic idea consists in using an explicit multistep method to compute a better initial guess and take advantage of this value by employing an implicit multistep method within a fixed-point scheme. The procedure is then divided into two parts: the first one is the prediction phase, where an explicit algorithm, named \textit{Predictor}, furnishes an adequate initial guess; the second one is the correction phase, where an implicit algorithm is invoked, possibly also several times, to realize the fixed-point scheme. The implicit method used is defined the \textit{Corrector} because acts on the predicted initial value. However, it is important to notice that the overall strategy is totally explicit because the predicted value is employed within the implicit method where the dependence on the incoming time instant appears. We will recall in details these concepts in Section~\ref{sec:algorithm}. The contents are organized as follows. 

In Section~\ref{sec:BK} we introduce the model and describe a particular version among those available as developments of the original one proposed by Burridge and Knopoff, mentioning the important connection with the Gutenberg-Richter law. In Section~\ref{sec:algorithm} we present the numerical algorithm and analyse the computational strategy used: we comment on the results of simulations by starting from the simplest case in order to increase the complexity and consider more articulated configurations. In Section~\ref{sec:increasing} we perform investigations of possible continuum limits of the essentially discrete Burridge-Knopoff model. Finally, in Section~\ref{sec:disc} a summary of our results is provided and perspectives on future work are discussed.
 
\section{The Burridge-Knopoff model}\label{sec:BK}
\selectlanguage{english}
 
\begin{figure}[h]
\centering
{\includegraphics[width=.60\columnwidth]{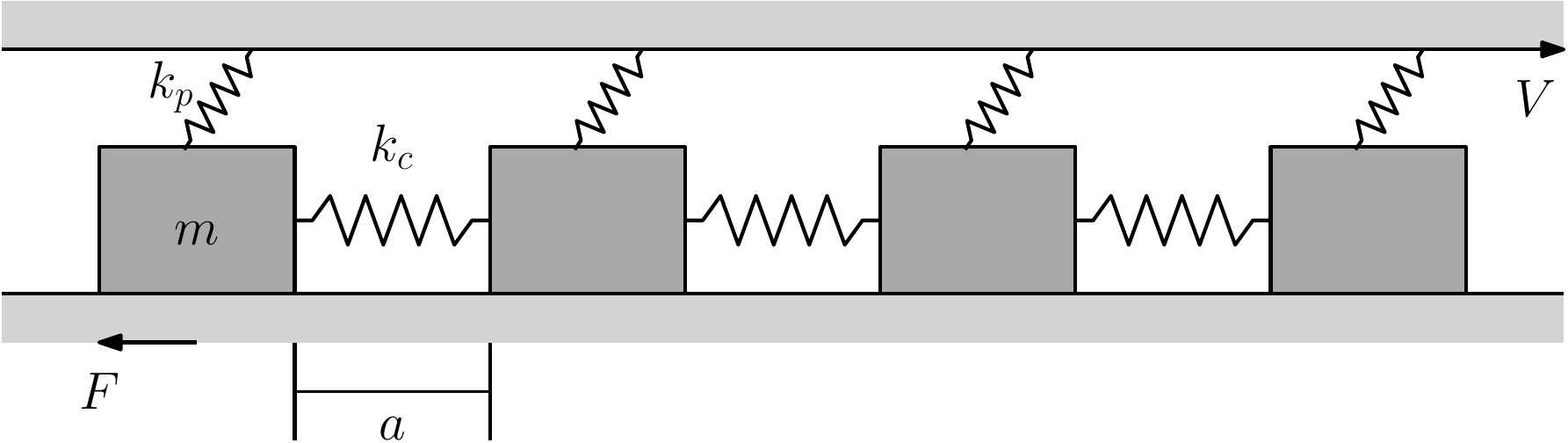}}\quad
\caption{Scheme of the Burridge-Knopoff model.}
\label{fig:cinqueblocchi}
\end{figure}

The system studied by Burridge and Knopoff is a spring-block model. Their purpose consists in trying to reproduce the typical dynamics which take place along an active fault. The goal is pursued through a discrete representation given by a chain of $N$ identical blocks, with mass $m$, mutually connected by linear springs with elastic constant $k_c$. A sort of one-dimensional array is generated (Fig.~\ref{fig:cinqueblocchi}). It is also possible studying the dynamics produced by a grid of blocks, within a multidimensional version of the system, thus focusing on a two-dimensional array \cite{MoriKawa2008twodimA,MoriKawa2008twodimB}.  
The blocks are supposed to rest on a rough surface, where $F$ is the friction, and connect to a moving upper plate by linear springs with elastic constant $k_p$. As regards the approximation of a real fault, the opposite sides of two different tectonic plates are assumed to be represented by the rough surface and the chain of blocks. The upper surface is supposed to be in motion, precisely at constant velocity $V$: this contribution induces a solicitation explainable thinking about the role of the external forces acting on a fault. It is assumed that the blocks are initially equally spaced and that the reciprocal distance is $a$. This means that $a$ does not explicitly appear within the equation of motion for the block $i$, which is 
\begin{equation}
\label{edoi}
m\ddot{x}_{i}=k_{c}(x_{i+1}-2x_{i}+x_{i-1})+k_{p}(Vt-x_{i})-F(\dot{x}_{i}),
\end{equation} 
where $x_{i}$ is the displacement from the initial equilibrium position. Let us investigate the structure of~\eqref{edoi} by analysing each contribution.
\subsection*{Internal elastic energy} As concerns the horizontal springs, it is assumed that a linear interaction takes place among the blocks. This is the conventional adjustment adopted within the Burridge-Knopoff model, but a linear coupling is not the only possibility. For instance, the eventuality of a nonlinear coupling is considered in \cite{Comte2002}. Due to the chain structure, producing two neighbors for each block, the internal elastic solicitation consists of two contributions, obviously except for the  masses at the edges (in this case adequate boundary conditions are required as we will discuss in Section~\ref{sec:algorithm}). By considering the elastic forces and recalling that the expression $x_{i}$ is associated to the displacements from equilibrium position, the contribution provided by springs with stiffness $k_c$ takes the form of the one-dimensional discrete Laplace operator.
\subsection*{External forces} We said above that the action of the external forces is realized within the model by the upper surface, in motion with constant velocity $V$. The blocks deal with this external element through the springs with stiffness $k_p$. So each mass is affected by another elastic solicitation besides that produced by the horizontal coupling. Of course, to quantify the vertical elastic force, it is necessary taking into account the elongation of springs caused by the upper plate. This consideration simply justifies the product $Vt$. It just has to combine this quantity with $x_{i}$ and $k_p$ according to the linear elasticity as in~\eqref{edoi}.
\subsection*{Friction} 
\begin{figure}[!h]
\centering
{\includegraphics[width=.60\columnwidth]{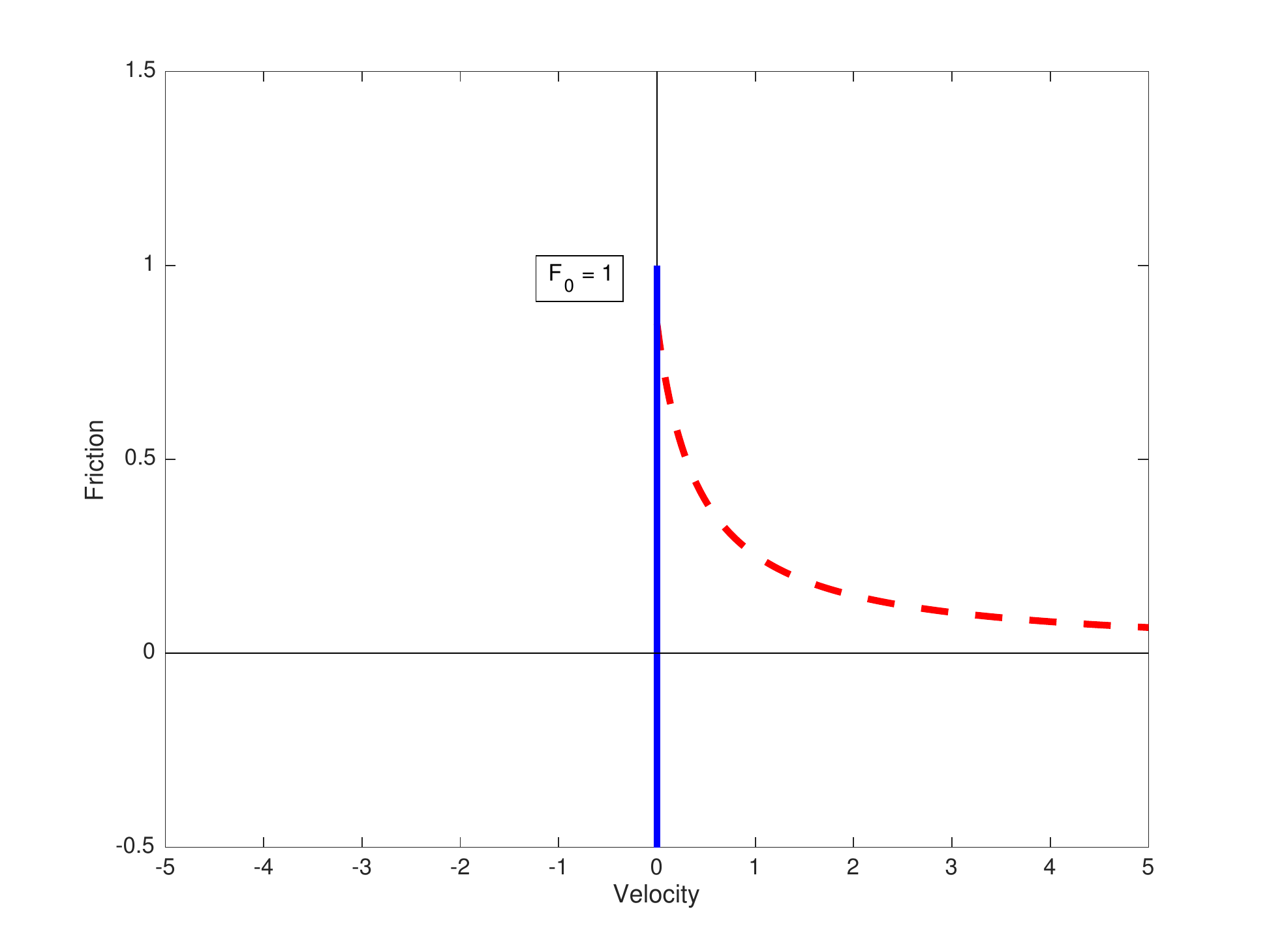}}\quad
\caption{The form of the friction law: the solid line is referred to the sticking friction, the dash one, to the slipping friction. $F_{0}$ is assumed to be unity.}
\label{fig:friction}
\end{figure}

The friction force $F(\dot{x}_{i})$ is velocity-dependent. This law allows to reproduce the typical stick-slip behaviour and introduces an essential instability inside the model. It is possible to distinguish different forms of friction, For instance, the Dieterich-Ruina friction law \cite{Erick2008,Erick2011,Ruina1983,Scholz1998}, the Coulomb friction law used by Muratov in \cite{Mura99} or the velocity-weakening friction proposed by Carlson and Langer \cite{CarlLang89properties,CarlLang89,CarlLang91intrinsic}. We adopt this last point of view in this paper. It is important to point out that another choice can be made between two different qualitative behaviours, the so-called asymmetric and symmetric versions, whose main difference is the constraint of non-negative velocity assumed in the asymmetric version. This means that back slip is inhibited for each block. We assume this constraint according to \cite{CarlLang91intrinsic,Langer96,MoriKawa2006,MoriKawa2008,Saito2007,Xia2008} and adopt the following (multi-valued) functional form
\[
F(v)=
\begin{cases}
\displaystyle \frac{F_{0}(1-\sigma)}{1+\frac{2\alpha v}{1-\sigma}} & \text{if $v>0$}  \\
(-\infty,F_{0}] & \text{if $v \leq 0$,} 
\end{cases}
\] 
where $v=\dot{x}_{i}$.
This double structure is easily understandable because a law based on the stick-slip dynamics must exhibit a discontinuity, as a consequence of the alternation between sticking and sliding motion for each block. The back sliding motion is forbidden by imposing $F(v)=-\infty $ for $v<0$. The value $F_{0}$ corresponds to the maximum static frictional force, so the static friction formally may range in the interval $(-\infty,F_{0}]$. During a sticking period the elastic resultant force acting on a block is perfectly balanced by the static friction, which means no motion. When the resultant force exceeds the threshold $F_{0}$, a slipping period starts with dynamic friction. Friction becomes weaker now, it decays monotonically to zero, as the velocity increases (see Fig.~\ref{fig:friction}). 
Another important feature of the friction law is the role of the parameters $\sigma$ and $\alpha$. The first one quantifies a small drop of the friction at the end of a sticking period; the second one provides informations about the decreasing of the dynamic friction force in relation to the increasing of the sliding velocity. 

One of the most interesting features of the Burridge-Knopoff model consists in the reproducibility of some important properties related to complex phenomena as real earthquakes, although the system exhibits a relatively simple structure. Among these typical behaviours, the {\it Gutenberg-Richter law} plays a significant role and can be used as a powerful instrument to assess the reliability of the model. This power law establishes that, in a seismic zone, the relationship between the number $N$ of earthquakes with intensity greater than or equal to a given magnitude $M$ and the magnitude itself has the form
\begin{equation}
\label{G.R.}
\log_{10}{N}= \bar{a}-bM.
\end{equation}
for some parameters $\bar{a}$ and $b$.
In order to represent the rate of seismic events, by introducing the total number of events expressed as $N_{T}=10^{\bar{a}}$, 
it is possible to reformulate the relationship~\eqref{G.R.} as
\begin{equation}
\label{G.R.bis}
\log_{10}({N/N_{T}})= -bM.
\end{equation}
This substitution allows us to understand the meaning of the quantity\footnote{This value is often called $a$: in this paper we adopt the notation $\bar{a}$ to avoid ambiguity with the quantity used to indicate the distance among the blocks in the Burridge-Knopoff model.} $\bar{a}$ in terms of total seismicity rate of an active zone. 
Finally, as regards the parameter $b$, in real situations its values are usually very close to $1$ in seismic zones \cite{deArcangelis2016}. 

It would be very interesting to make comparisons also within the aftershocks field, by studying the {\it Omori law} \cite{deArcangelis2016}. However, at least in a such simple version of the Burridge-Knopoff model as the current one, it is impossible recognizing aftershocks sequences, as pointed out in \cite{Kawamura2012}. Further contributions, as viscosity, would be required.

\section{Numerical algorithm and its reliability}\label{sec:algorithm}
\selectlanguage{english}
In order to justify our simulations, before showing the results, we discuss the computational strategy adopted and the numerical algorithm chosen. 

Once the number of blocks $N$ is established, it is possible to obtain a differential system composed by equations like \eqref{edoi}, namely
\[
\begin{sistema}
\dot{x}_{i}= y_{i} \\
m\dot{y}_{i}= k_{c}(x_{i+1}-2x_{i}+x_{i-1})+k_{p}(Vt-x_{i})-F(y_{i}),
\end{sistema}
\]
for $i=1$, \ldots, $N$, where we assume that $x_{0}=x_{1}$ and $x_{N+1}=x_{N}$ for the boundary conditions \cite{BurrKnop67}. 
\subsection*{Initial conditions} In order to avoid a periodic evolution, with the aim of reproducing realistic local tension along a fault, we assign small random displacements from the equilibrium positions for each mass as in \cite{XionKikuYama15}. The blocks are supposed to be at rest, so we put zero velocities. Remembering that the blocks are initially equally spaced with distance $a$, the equilibrium positions are $P_{i}=a(i-1)$ for $i=1, \dots, N$. Because zero velocities are imposed, all blocks are initially stuck. However, if a simulation would have started with the actual initial conditions, some irregular dynamical motion would be recognized, due to the action of spring force and friction force. On the contrary, we wish to appreciate a realistic charge cycle. That is why we identify the next incoming time of global stick, $\bar{t}$, and select this one as initial time. This implies that the original initial conditions, and corresponding perturbations, must be updated in $\bar{t}$: the simulation is now ready to be restarted by setting $t=0$. 

\subsection*{Stick-slip detection} To identify a time of stick for a block within the numerical code, we use a criterion based on both the resultant force and velocities: a block is  stuck if and only if the elastic resultant force is less than the maximum static frictional force and the velocity is equal to zero. Obviously, it is very difficult detecting an absolutely zero velocity in simulations so that we use a workaround: because back slip is inhibited, sign changes of the velocity are interpreted as the tendency of being stuck, so negative values are suppressed and replaced by zero values. That is why we do not need to introduce a threshold parameter to create a range for the zero value as it is often done working with the asymmetric friction law. 

\subsection*{Seismic events} Talking about the statistical properties of the model, specifically referred to the Gutenberg-Richter law, we have to assume an operative definition to judge whether a seismic event is happening: an earthquake occurs when a blocks starts to slip and ends only where all the blocks are stuck again. This definition implies that, during an event, a block can slip and become stuck alternately; moreover, the elastic coupling produces a sort of propagation along the chain of masses, because a block can trigger the slipping of its neighbors. In order to quantify the magnitude of an earthquake, we introduce the following definition for the magnitude $M$:
\begin{equation}
\label{magni}
M=\log_{10}\biggl(\sum_{i=1}^{N}\Delta x_{i}\biggr),
\end{equation}
where $\Delta x_{i}$ is the cumulative displacement of the block $i$ during a given earthquake.
\subsection*{Numerical adjustment} As regards the numerical integration of the Burridge-Knopoff model with velocity-weakening friction, various methodologies have been employed by using either explicit methods, such as explicit Runge-Kutta \cite{MoriKawa2006,MoriKawa2008twodimA,MoriKawa2008twodimB,
MoriKawa2008,Xia2008}, or implicit methods, such as Implicit Euler \cite{XionKikuYama15}. 

We adopt a Predictor-Corrector strategy. Let us explore this procedure \cite{Quarteroni}. First of all we start by considering a general implicit multistep method, For instance, by selecting the Adams-Moulton family from the Adams methods. The following equation groups all the Adams methods,
\begin{equation}
\label{Adams}
y_{n+1}=y_{n}+h\sum_{j=-1}^{p}b_{j}f_{n-j}.
\end{equation}
If $b_{-1}\neq 0$, an implicit method, named Adams-Moulton, is generated; otherwise, when  $b_{-1}=0$, an explicit method, named Adams-Bashforth, is obtained. That is why we assume $b_{-1}\neq 0$. In~\eqref{Adams} $y_{n}$ indicates the approximate solution evaluated at the time $t_{n}$; $f_{n-j}$, more explicitly $f(t_{n-j},y_{n-j})$, is the vector field; $h$ is the step size; $b_{j}\in\mathbb{R}$; $p\in\mathbb{N}$ is used to quantify the number of steps of the method, precisely $p+1$, without including the implicit part associated to $j=-1$. We recall that the Adams methods are derived from the integral representation of the Cauchy's problem for a given differential system, namely
\begin{equation*}
\label{IntegralCauchy}
x(t)=x_{o}+\int_{t_{0}}^{t} f(s,x(s))\, ds,
\end{equation*}
by using interpolating polynomials in the Lagrange form to approximate the vector field~\cite{Quarteroni}. In order to solve a Cauchy's problem by using an implicit method such as~\eqref{Adams} it is necessary to approach a nonlinear equation. We can rewrite~\eqref{Adams} as follows 
\begin{equation}
\label{fixedpoint}
y_{n+1}=y_{n}+h\sum_{j=-1}^{p}b_{j}f_{n-j}=\Phi(y_{n+1}).
\end{equation}  
By taking advantage of~\eqref{fixedpoint} we can adopt fixed-point iterations and thus solve the nonlinear equation. For $k=0,1\dots,$ we get
\begin{equation}
\label{fixedpointbis}
y^{(k+1)}_{n+1}=\Phi(y^{(k)}_{n+1}).
\end{equation}
However the procedure triggered by~\eqref{fixedpointbis} requires several function evaluations to achieve convergence, due to the iterations needed. The idea behind a Predictor-Corrector strategy, which is inspired by the purpose of reducing the computational cost, is to compute a good initial guess for the fixed-point iterations by recalling an explicit multistep method. This method, called \textit{Predictor}, provides an adequate guess to be used within the fixed-point scheme~\eqref{fixedpointbis} generated by the implicit algorithm~\eqref{Adams}. The implicit method, named \textit{Corrector}, can be invoked $m$ times, with $m\geq1$. When $m>1$ the procedure is called \textit{Predictor-Multicorrector}. In this paper we choose $m=1$, so we will continue using simply the wording Predictor-Corrector. The algorithm produced by starting from the Adams methods can be summed up as follows
\[
\begin{sistema}
\begin{aligned}
\textit{Predict}&: &y^{(0)}_{n+1}&=y^{(1)}_{n}+h\sum_{j=0}^{\bar{p}}\bar b_{j}f^{(0)}_{n-j} \\
\textit{Evaluate}&: & f^{(0)}_{n+1} &= f(t_{n+1},y^{(0)}_{n+1}) \\
\textit{Correct}&: &y^{(1)}_{n+1}&=y^{(1)}_{n}+hb_{-1}f^{(0)}_{n+1}+h\sum_{j=0}^{p}b_{j}f^{(0)}_{n-j},
\end{aligned}

\end{sistema}
\]
where the \textit{Evaluation step} of the vector field $f$ is included. The superscript $(0)$ denotes the guess provided by the Predictor, the superscript $(1)$, instead, indicates the values furnished by the Corrector. The abbreviation usually employed for the overall procedure is \textit{PEC}. We notice that a Predictor-Corrector strategy, also in the general case $m\geq1$, is by construction totally explicit. As regards our numerical adjustment, we adopt the Predictor-Corrector technique in a bit different form, called \textit{PECE}, in which a further evaluation of $f$ is performed at the end of the sequence. Moreover, the second-order Adams-Bashforth scheme (AB2) is used as Predictor, while the third-order Adams-Moulton method (AM3) is chosen as Corrector. We thus obtain
\[
\begin{sistema}
\begin{aligned}
\textit{Predict} &: & y^{(0)}_{n+1} &= y^{(1)}_{n} +\frac{h}{2}[3f^{(1)}_{n}-f^{(1)}_{n-1}]  \\
\textit{Evaluate}&: & f^{(0)}_{n+1} &= f(t_{n+1},y^{(0)}_{n+1}) \\
\textit{Correct} &: & y^{(1)}_{n+1}      &= y^{(1)}_{n} +\frac{h}{12}[5f^{(0)}_{n+1}+8f^{(1)}_{n}-f^{(1)}_{n-1}] \\
\textit{Evaluate}&: & f^{(1)}_{n+1}      &= f(t_{n+1},y^{(1)}_{n+1}).

\end{aligned}

\end{sistema}
\]
We point out that the order, $q$, of the \textit{PECE} procedure, can be computed as follows
\begin{equation*}
\label{order}
q=\min(q_{p}+1,q_{c}),
\end{equation*}
where $q_{p}$ and $q_{c}$ are the orders of the Predictor and the Corrector steps, respectively. 
Therefore, we have generated an overall third-order method. 

Evaluating the pros and the cons of the Predictor-Corrector technique, the main advantages consist, firstly, in avoiding to solve an implicit system at every step, whose size would increase with the number of blocks, and, secondly, in ensuring a stronger stability when compared to standard explicit methods. On the other hand, we have to adopt a small time-step to achieve a good approximation of the solution to the Burridge-Knopoff model. According to \cite{MoriKawa2008,Xia2008}, we choose $h=0.001$ as constant step.

\subsection*{One block}
\begin{figure}[h]
\centering
{\includegraphics[width=.60\columnwidth]{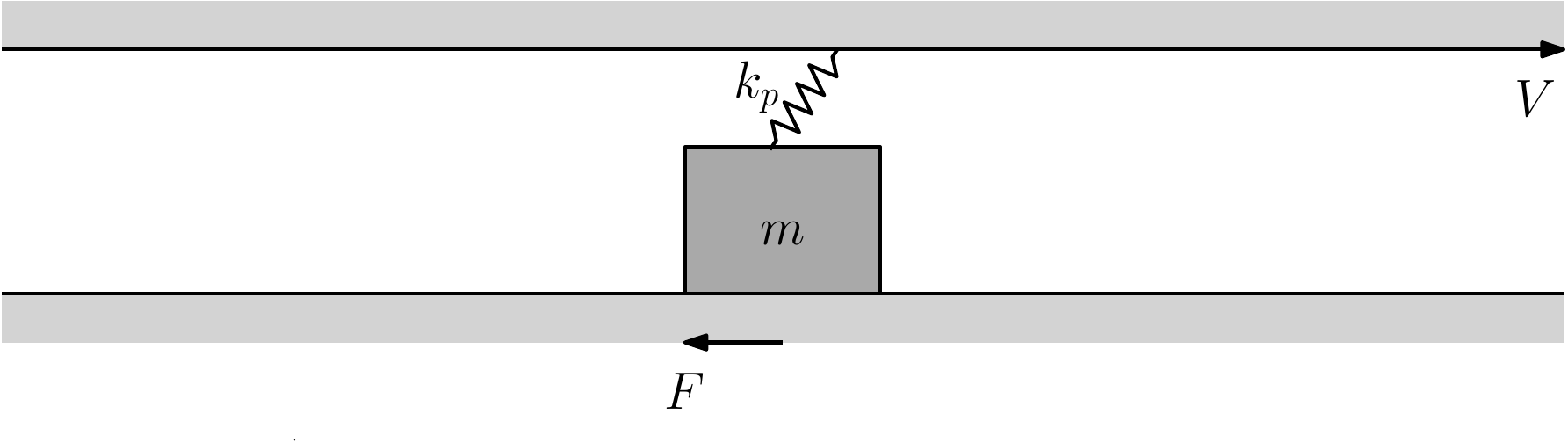}}\quad
\caption{The system involving a single block.}
\label{fig:single}
\end{figure}

Let us get into the simulations now.
As the simplest case, we examine the evolution of the system in which only one block is involved (see Fig.~\ref{fig:single}). 
The aim is to become familiar with the specific trend of the stick-slip dynamics. The equation of motion can be easily deduced from~\eqref{edoi} by omitting the elastic term associated to the horizontal connecting springs, because in this configuration adjacent blocks are not included, so that we get
\begin{equation}
\label{OneBlock}
m\ddot{x}=k_{p}(Vt-x)-F(\dot{x}).
\end{equation} 
For the values of the parameters involved in the model, we follow the work by Saito and Matsukawa \cite{Saito2007} and adopt the list shown in Table~\ref{tab:param}. We deduce from this table all the values useful to integrate~\eqref{OneBlock}. Finally, for the remaining parameters, we put $\sigma=0.01$ according to \cite{MoriKawa2006} and arbitrarily choose $\alpha=1$. In the follow-up we will discuss carefully the role of the quantity $\alpha$, which is very significant within the configurations involving lots of blocks.

As initial conditions we simply impose $x(0)=0$ and $\dot{x}(0)=0$ without adopting artifices as those mentioned previously, very useful in the case of more blocks. It is obviously possible to assume a more realistic small displacement for $x(t)$ in $t=0$ but the evolution would not change its qualitative behaviour, the only difference consisting in the duration of the first stick period. When only one block is involved, indeed, the motion exhibits a periodic trend. In order to avoid this kind of dynamics, it is crucial introducing more blocks within the system.
\begin{table}[!h]
\caption{Values of quantities involved in the simulations with one block.}
\label{tab:param}
\centering
\begin{tabular}{*{5}{c}}
\toprule
\multicolumn{5}{c}%
{\textbf{Parameters}}  \\   
\midrule
$\mathbf{m}$ &$\mathbf{k_{p}}$ &$\mathbf{k_{c}}$ &$\mathbf{V}$ &$\mathbf{F_{0}}$  \\  
\midrule  
  $1$      & $1$       & $60$      & $0.001$ & $1$   \\
\bottomrule
\end{tabular}
\end{table}

\begin{figure}[!h]
\centering
{\includegraphics[width=.60\columnwidth]{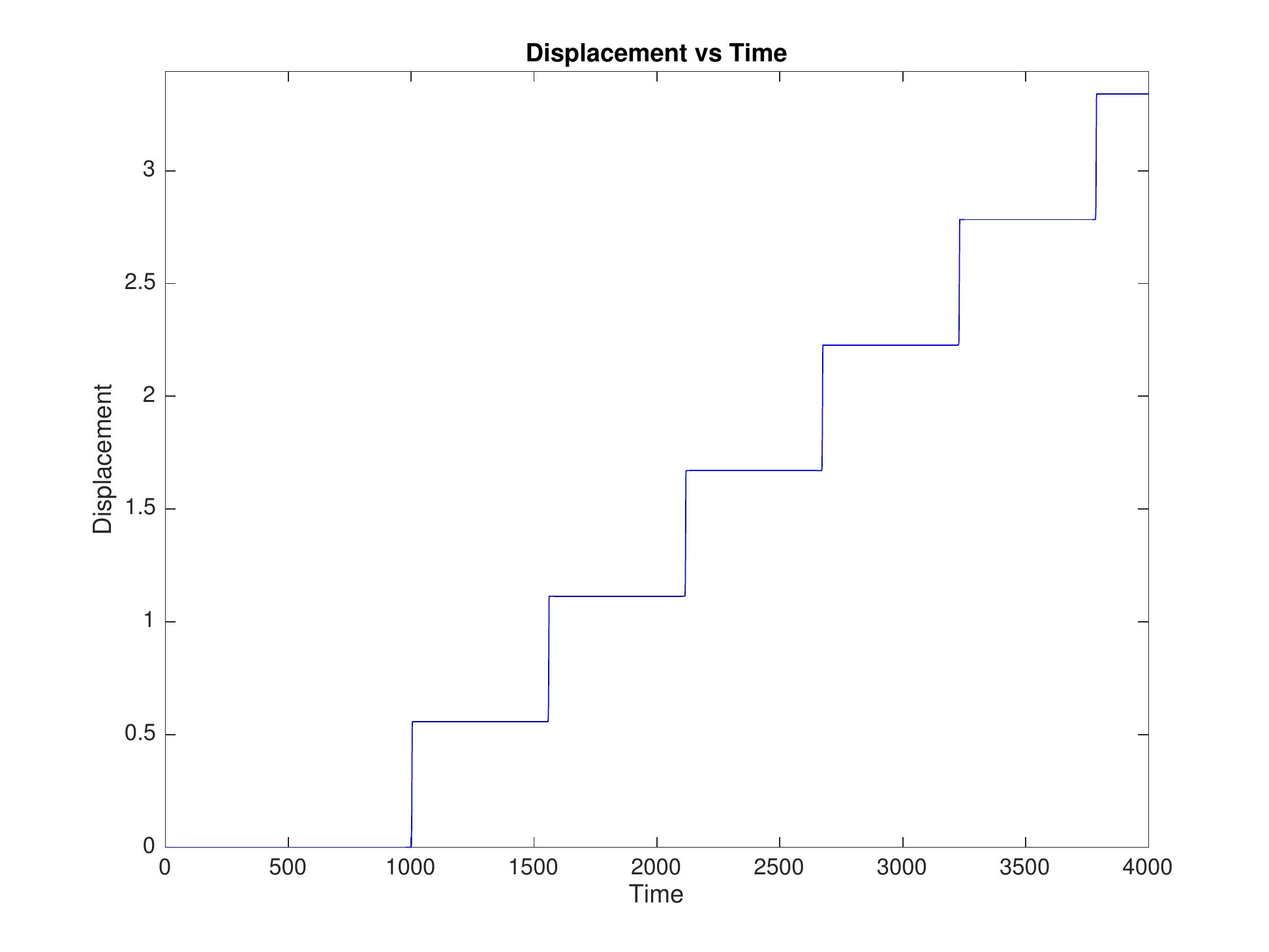}}\quad
\caption{Displacement for a single block. All the values used to perform the simulation are available in Table~\ref{tab:param}.}
\label{fig:SingleDisplacement}
\end{figure}

\begin{figure} [!h]
\centering
\subfloat[][\emph{Velocity graph.}]
{\includegraphics[width=.487\columnwidth]{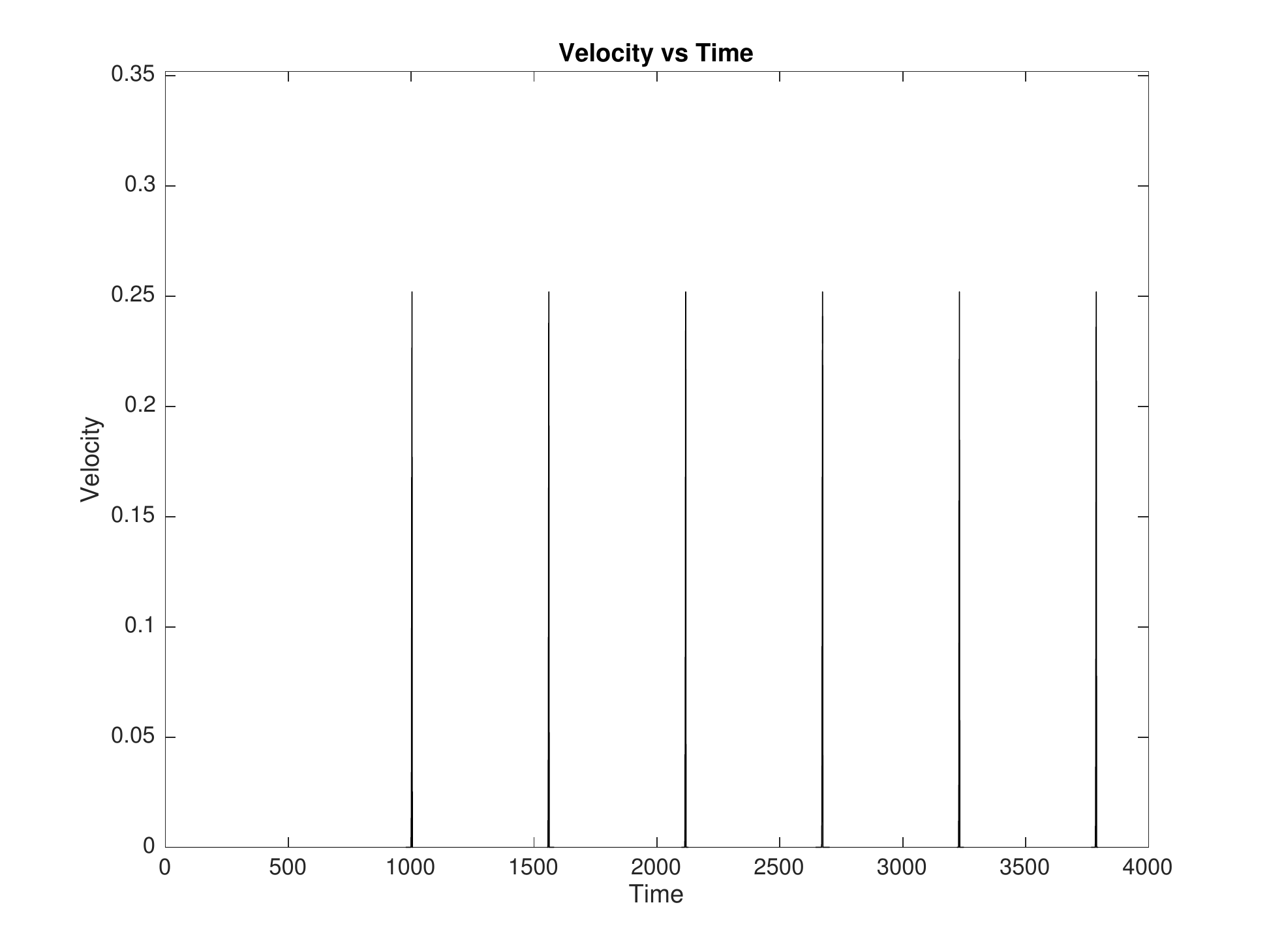}}\quad
\subfloat[][\emph{Phase portrait.}]
{\includegraphics[width=.487\columnwidth]{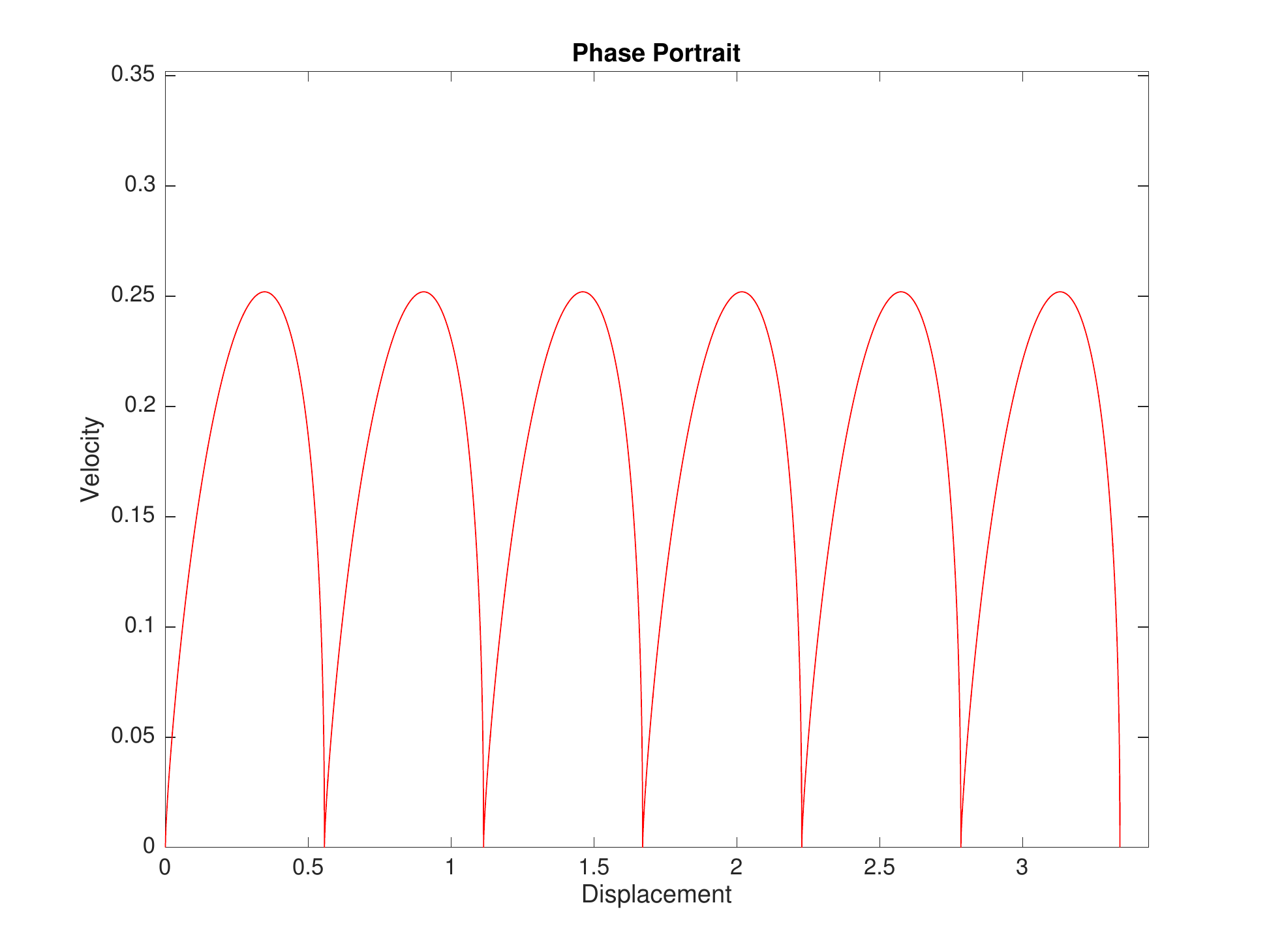}}\quad
\caption{Velocity and phase portrait for a single block. All the values used to perform the simulation are available in Table~\ref{tab:param}.}
\label{fig:VelocityPhase}
\end{figure}

Figs.~\ref{fig:SingleDisplacement} and~\ref{fig:VelocityPhase} show the results in the case of a single block. In absence of adjacent masses, the motion tends to be periodic. We recognize the alternation between sticking and slipping periods from the qualitative behaviour of the graph in Fig.~\ref{fig:SingleDisplacement}: a steep trend characterizes the sliding motion, in opposition to the flat one produced when there is not motion. 
When the block is sliding, its velocity achieves some pronounced peaks as shown in Fig.~\ref{fig:VelocityPhase}a. In Fig.~\ref{fig:VelocityPhase}b a qualitative summary is provided by the phase portrait. Although this kind of system is very far from being an accurate representation for seismic events, because lots of other contributions would be required, on a geophysical level we could think about a slipping period as an earthquake and a sticking one as a charge cycle.
\subsection*{A bit more complex configuration} 
\begin{figure} [!h]
\centering
\subfloat[][\emph{Displacements as a function of time.}]
{\includegraphics[width=.487\columnwidth]{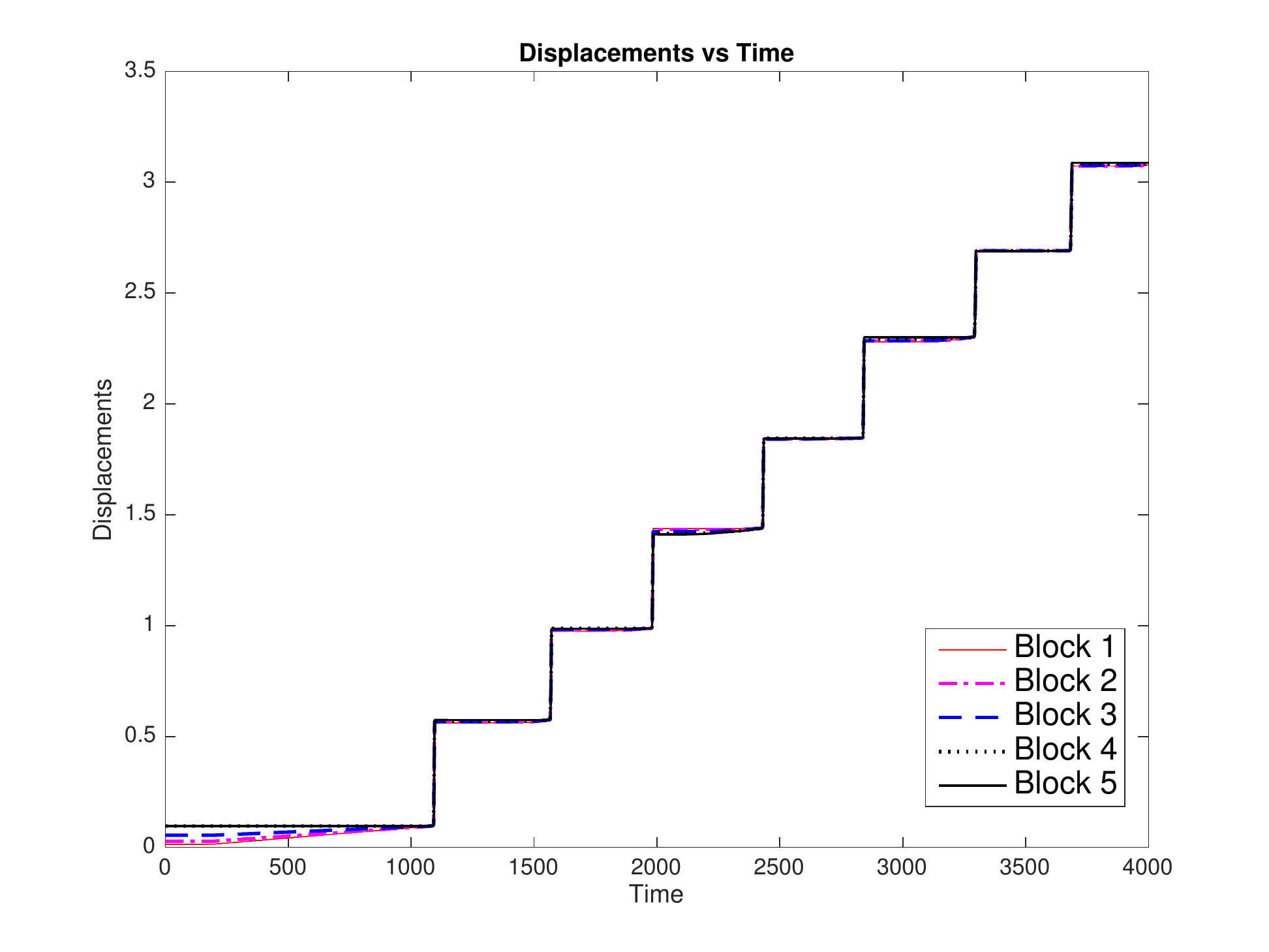}}\quad
\subfloat[][\emph{Zoom-in of the previous graph.}]
{\includegraphics[width=.487\columnwidth]{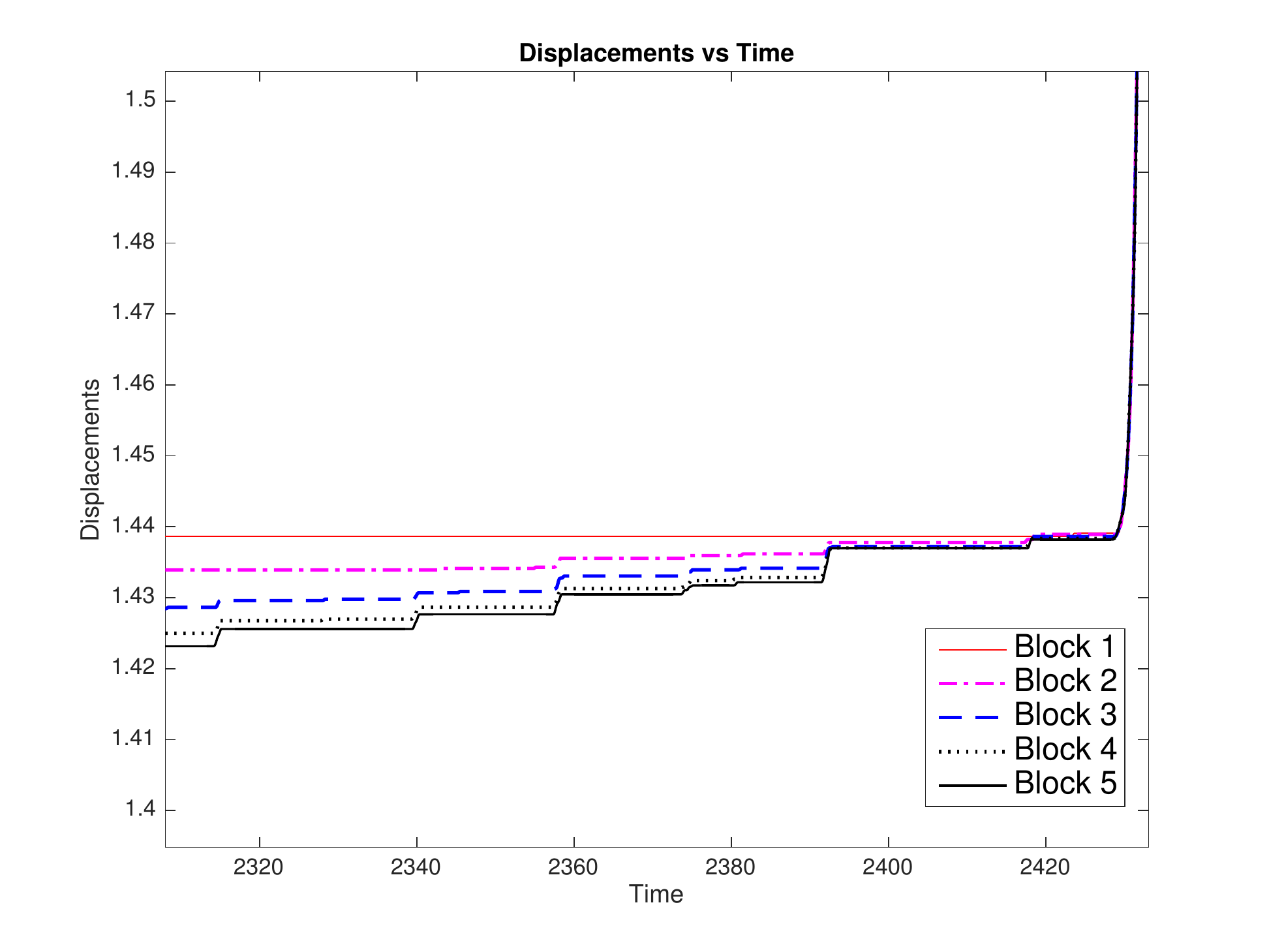}}\quad
\caption{A system involving five blocks. The graph of the displacements as a function of time
and its zoom-in for the parameters listed in Table~\ref{tab:param}.}
\label{fig:Five}
\end{figure}

\begin{figure} [!h]
\centering
\subfloat[][\emph{Velocity for the first block.}]
{\includegraphics[width=.487\columnwidth]{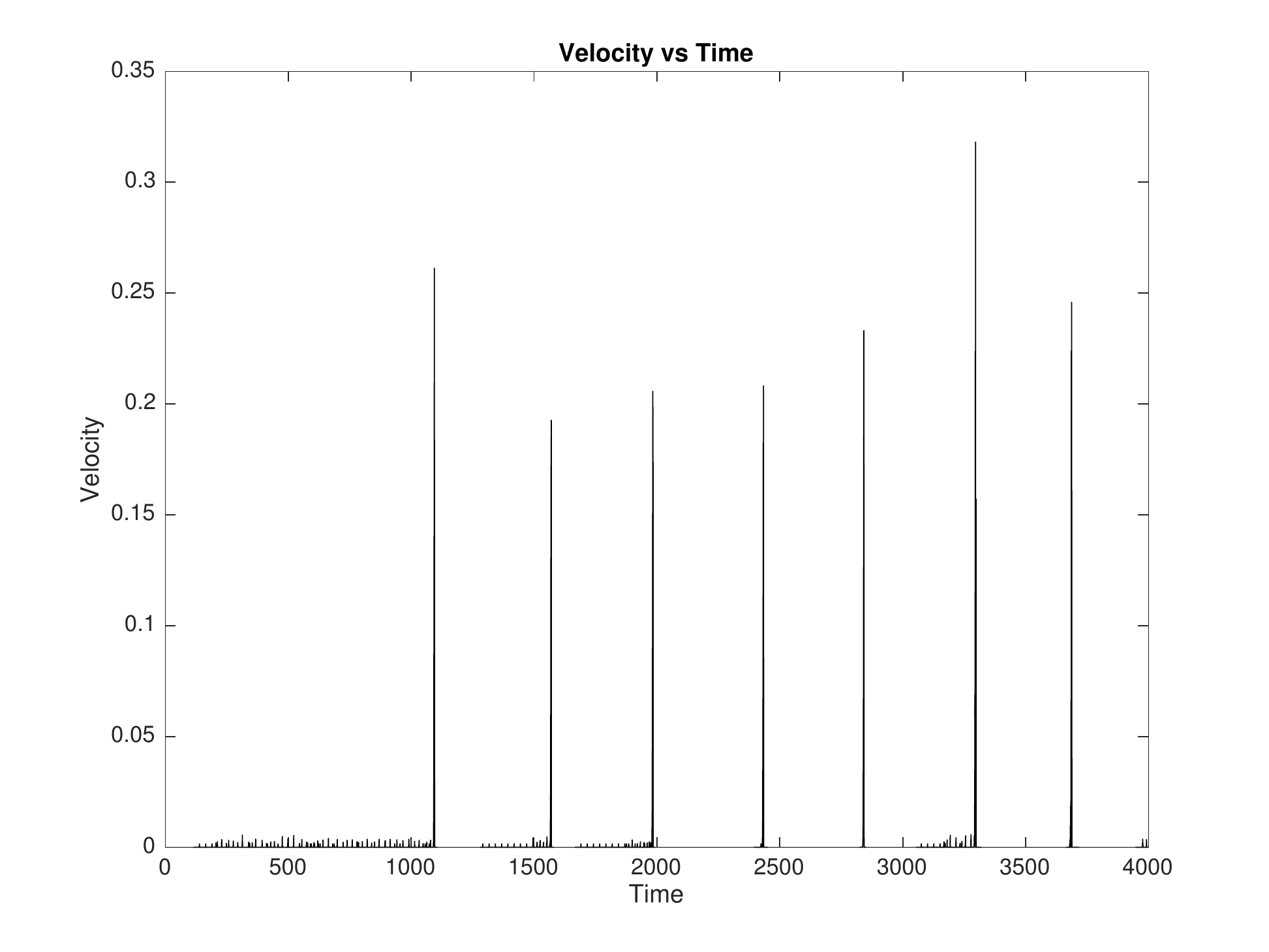}}\quad
\subfloat[][\emph{Phase portrait for the first block.}]
{\includegraphics[width=.487\columnwidth]{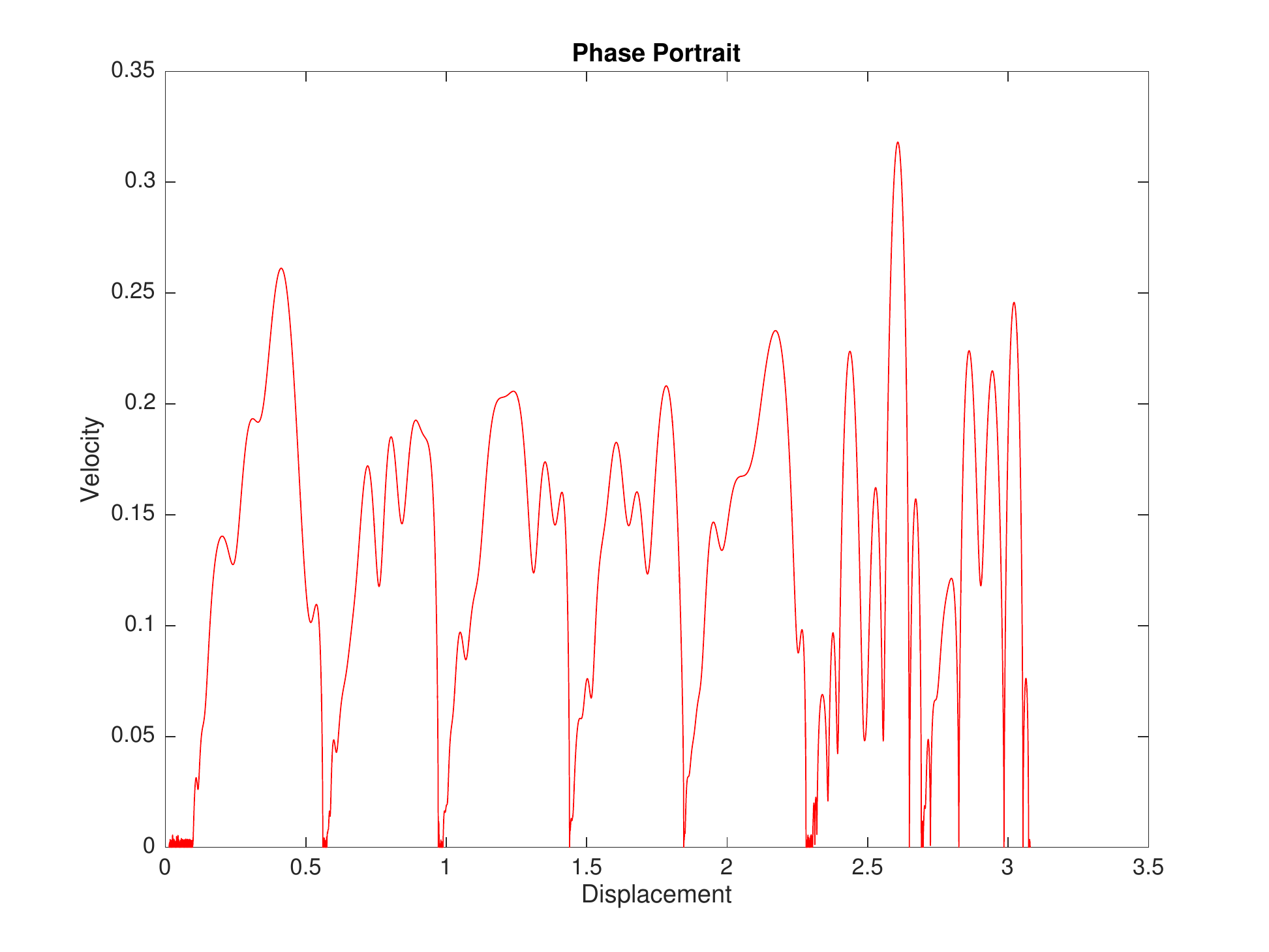}}\quad
\caption{Block $1$ has been taken as sample for the simulation including five blocks. 
The graph of velocity as a function of time and corresponding phase portrait are similar for the remaining blocks.}
\label{fig:bl1}
\end{figure}

We proceed to examine a model involving five masses. By adopting the same parameters used in the case of a single block, listed in Table~\ref{tab:param}, and assuming as initial conditions random, small displacements updated at the incoming time of global stick, as described before, we simulate the evolution. The associated differential system is defined by using~\eqref{edoi} and paying attention to include the boundary conditions. For instance, the equation for the first block becomes
\begin{equation*}
m\ddot{x}_{1}= k_{c}(x_{2}-x_{1})+k_{p}(Vt-x_{1})-F(\dot{x_{1}}).
\end{equation*}
In Fig.~\ref{fig:Five}a we plot the displacements of the blocks, while in Fig.~\ref{fig:Five}b a zoom-in for these trajectories is provided. It is possible to recognize some great slipping phases in which all the blocks are involved; on the other hand, talking about the flatter trends, very small displacements happen. The pronounced peaks correspond to the most powerful shocks allowed for such a limited configuration. Finally, looking carefully at the zoom-in Fig.~\ref{fig:Five}b, a sequence of narrow events happening before one of the peaks mentioned above is captured: these small earthquakes can be interpreted as foreshocks.
Although five blocks are not enough to exhibit a satisfying dynamics, it has been helpful to investigate the results in a qualitative way. As can be observed, indeed, the behaviour is certainly more complex and nontrivial than in the case of a single block. In order to support this point of view, in Fig.~\ref{fig:bl1} we take one block as sample and plot the velocity and the phase portrait (trends are very similar for the remaining blocks). So more facets and details about the qualitative trend are pointed out.  
\subsection*{More blocks and the Gutenberg-Richter law} 
\begin{figure}[!h]
\centering
{\includegraphics[width=.60\columnwidth]{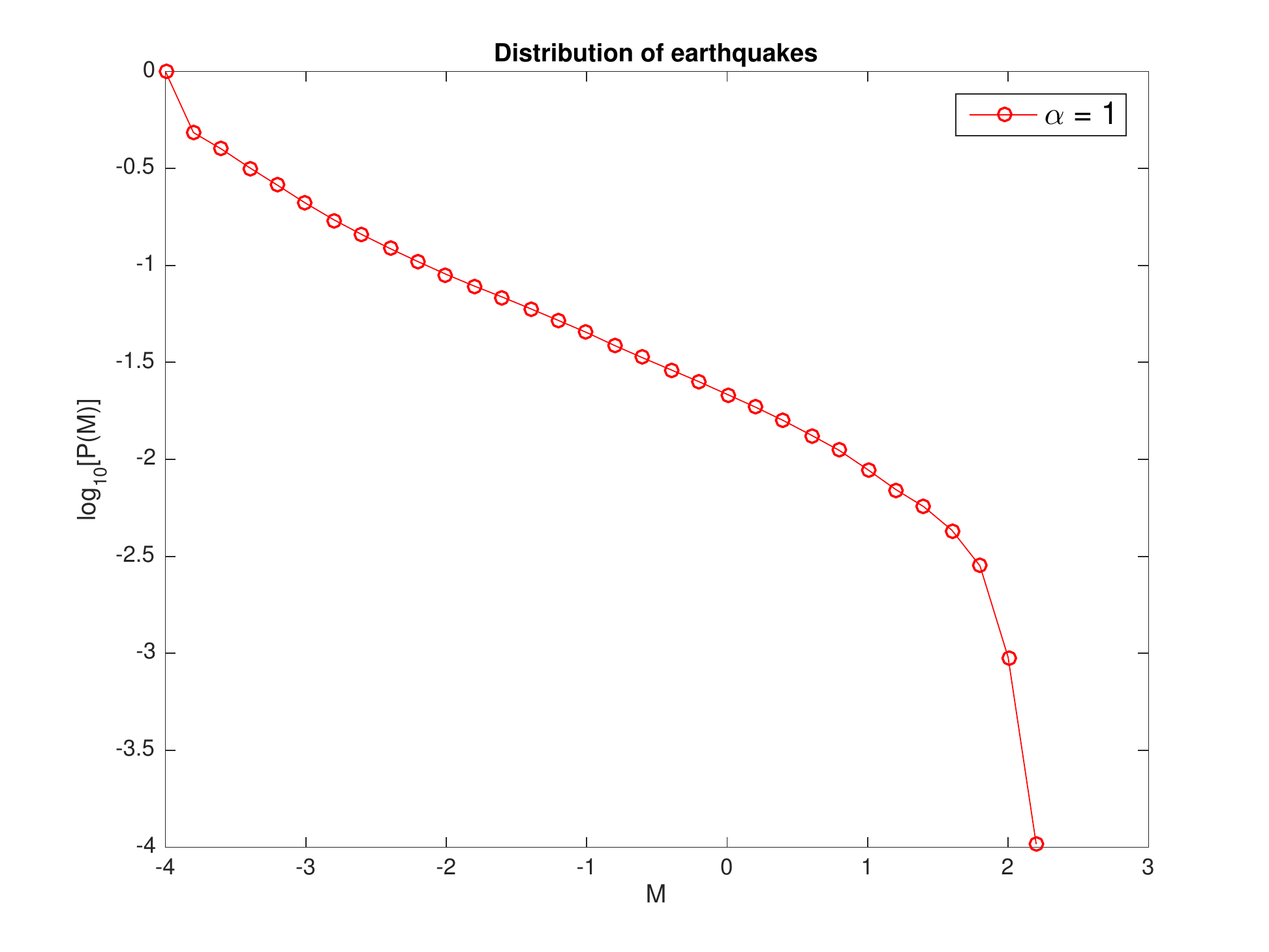}}\quad
\caption{Earthquakes distribution for a system with $200$ blocks: graph of the function $M\mapsto \log_{10}[P(M)]$ with $\alpha=1$.
The simulation is realized by using the parameters listed in Table~\ref{tab:param2}.}
\label{fig:Alfa1bl200GR}
\end{figure}

The next steps in our analysis are aimed at arguing the reliability of the Burridge-Knopoff model using the Gutenberg-Richter law exhibited in~\eqref{G.R.bis}. In order to achieve this purpose, we increase the number of blocks and consider a system including two hundred blocks, so $N=200$. In Table~\ref{tab:param2} the parameters used in this last part of Section~\ref{sec:algorithm} are listed.
\begin{table}[!h]
\caption{Values of quantities involved in the simulations with several blocks.}
\label{tab:param2}
\centering
\begin{tabular}{*{7}{c}}
\toprule
\multicolumn{7}{c}%
{\textbf{Parameters}}  \\   
\midrule
$\mathbf{m}$ &$\mathbf{k_{p}}$ &$\mathbf{k_{c}}$ &$\mathbf{V}$ &$\mathbf{F_{0}}$  &$\mathbf{\sigma}$  & $\mathbf{\alpha}$\\  
\midrule  
  $1$      & $1$       & $100$      & $0.001$ & $1$ & $0.01$ & $\{1,1.5,2,3,4\}$  \\
\bottomrule
\end{tabular}
\end{table}
We collect informations about the seismic events generated by using the criterion described before to distinguish when an earthquake is happening into the simulation. As regards the magnitude, the relationship~\eqref{magni} provides a quantitative definition. First of all, to discuss the results, it is absolutely necessary focusing on the role of the parameter $\alpha$ introduced by the friction law \cite{CarlLang89,Langer96,MoriKawa2006,Xia2008}. As said in Section~\ref{sec:BK}, $\alpha$ expresses the rate of slipping friction decreasing on increasing the sliding velocity. As a result, if $\alpha$ decreases, friction becomes more dissipative; on the contrary, larger values of $\alpha$ mean less dissipation because the slipping friction decreases more quickly with sliding velocity. On a quantitative level, the value $\alpha=1$ is an important threshold, since values of $\alpha$ less than unity preventing system from exhibiting great earthquakes, that is why $\alpha=1$ is a sort of lower bound.  Moreover, different behaviours are noticeable by assuming $\alpha=1$ and $\alpha>1$: let us investigate this point. By following \cite{XionKikuYama15} we introduce the earthquake distribution $P(M)$, which is the ratio between the number of earthquakes greater than or equal to a given magnitude $M$ (see~\eqref{magni}) and the total number of events $N_{T}$. Operatively we classify magnitudes by establishing the belonging to different ranges such as $[M, M+dM]$, where $dM$ is fixed to be equal to $0.2$. According to~\eqref{G.R.bis}, we represent the distribution of earthquakes by the graph of the function $M\mapsto \log_{10}[P(M)]$. In the first case considered $\alpha=1$ is employed and  the result is shown in Fig.~\ref{fig:Alfa1bl200GR}.

\begin{figure} [!h]
\centering
{\includegraphics[width=.60\columnwidth]{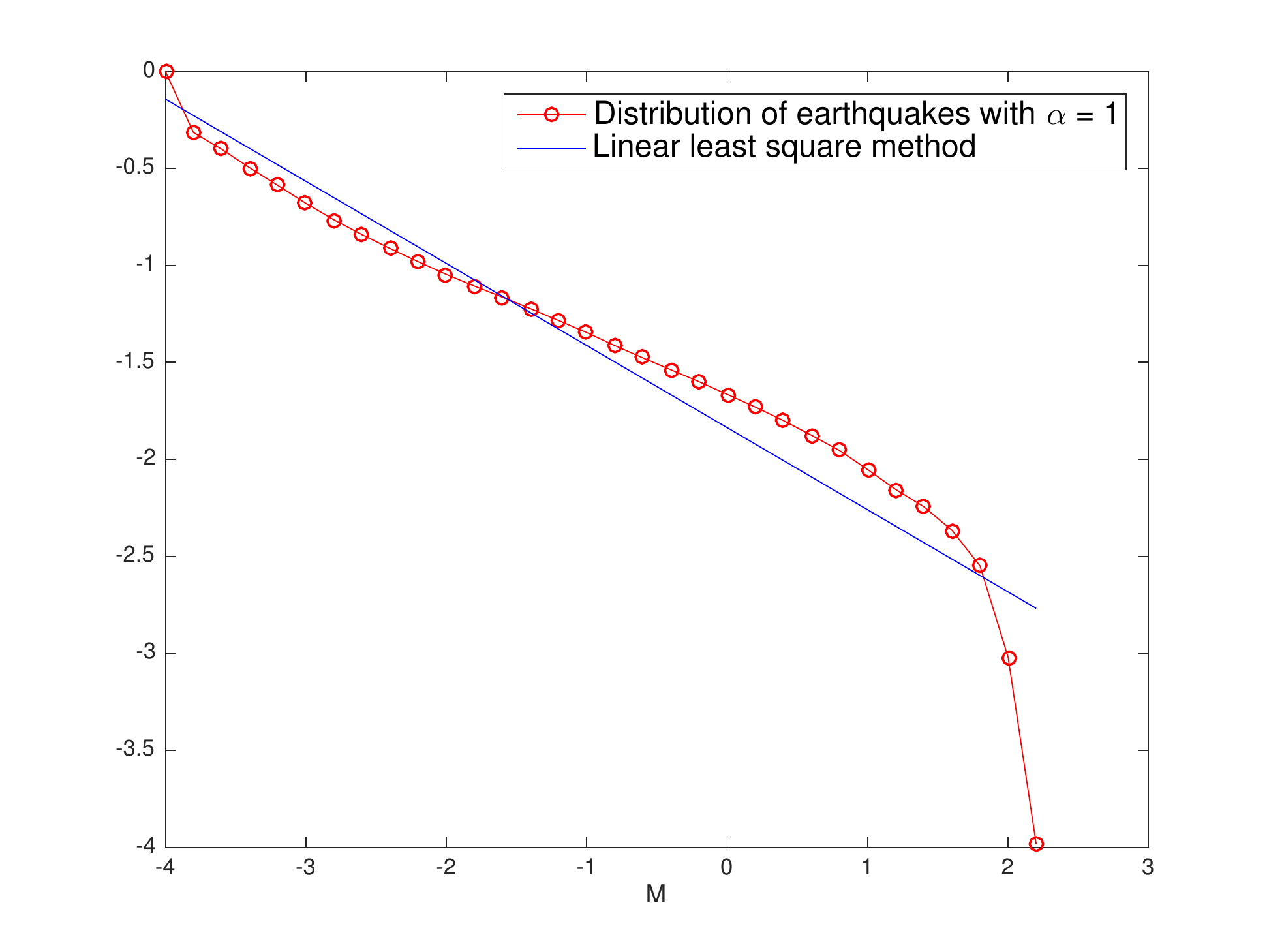}}\quad
\caption{Distribution of earthquakes and extrapolation by the linear least squares method when $\alpha=1$:
a trend very close to a straight line is recognized.}
\label{fig:MinQuad}
\end{figure}

By observing the graph, it is possible to recognize a behaviour very close to a straight line in the central part (between $M=-3.7$ and $M=1.7$, approximately), which can be interpreted coherently with the Gutenberg-Richter law. By adopting the linear least squares method (see Fig.~\ref{fig:MinQuad}) we estimate an exponent $B\simeq 0.42$ for the power-law trend. We point out that the exponent derived from the simulations of the Burridge-Knopoff model cannot be directly matched to the $b$-value appearing in the Gutenberg-Richter law~\eqref{G.R.bis}, indeed a rescaling would be required in order to make a comparison with the real data \cite{deArcangelis2016,Kawamura2012}.  We can take as a reference the relationship $b=\frac{3}{2}B$ described in \cite{Kawamura2012}. Moreover, according to the results in \cite{MoriKawa2006}, we notice that by varying the ratio between the stiffness of the springs, $k_{c}/k_{p}$, usually called $l^2$ in literature, a bit different exponent $B$ is computed: the lower the ratio the higher the exponent, of course preserving the constraint $k_{c}>k_{p}$. For instance, by imposing $l^2=36$ we find $B\simeq 0.45$ or choosing $l^2=9$ we obtain $B\simeq 0.54$. As regards the data providing the results plotted in Fig.~\ref{fig:MinQuad}, a ratio $l^2=100$ can be deduced from Table~\ref{tab:param2}. All these considerations make our $b$-value ranging in the interval $[0.63,0.81]$, that means a bit flatter slope for the graph: indeed it is less than the empiric value $b=1$ in~\eqref{G.R.bis}. Finding a flatter slope is consistent with other observations available in the literature, as in \cite{CarlLang91intrinsic,MoriKawa2006,XionKikuYama15}. 

At very small magnitude we notice a steep linear segment in the graph, probably caused by the discreteness of the model \cite{MoriKawa2006,XionKikuYama15}. 

\begin{figure}[!h]
\centering
{\includegraphics[width=.60\columnwidth]{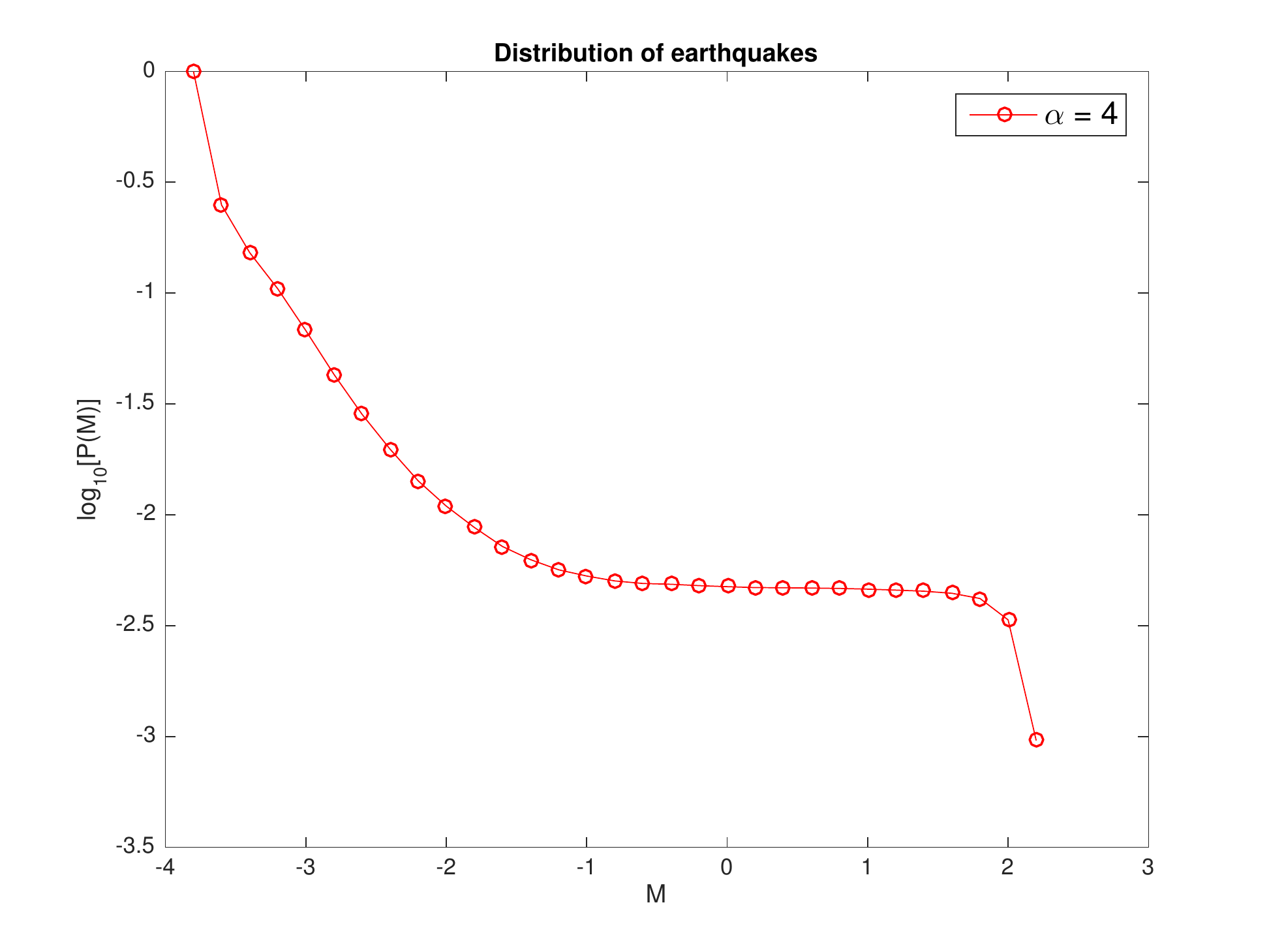}}\quad
\caption{Earthquakes distribution for a system with $200$ blocks: graph of the function $M\mapsto \log_{10}[P(M)]$ with $\alpha=4$.
All the other parameters employed are those in Table~\ref{tab:param2}.}
\label{fig:Alfa4bl200GR}
\end{figure}

By increasing $\alpha$, a different qualitative behaviour is provided by simulations. For instance, we assume now $\alpha=4$ (see Fig.~\ref{fig:Alfa4bl200GR}) in order to show the main differences, without neglecting to provide an accurate screening, by employing more values of $\alpha$, in what follows.  All the other parameters are the same as in the case $\alpha=1$. 

In Fig.~\ref{fig:Alfa4bl200GR} we recognize a deviation from the Gutenberg-Richter law at large magnitudes: it is noticeable a sort of peak structure; the trend close to a straight line persists instead in the middle-small range of magnitude, in agreement with the empirical expectation. Finally, at smallest magnitude, as in the case analysed before, a steep linear segment is observable due to the discreteness of the model. All these qualitative behaviours are consistent with previous works, for example \cite{CarlLang91intrinsic,MoriKawa2006,Saito2007,XionKikuYama15}. 

\begin{figure} [!h]
\centering
\subfloat[][\emph{$\alpha=1$}.]
{\includegraphics[width=.487\columnwidth]{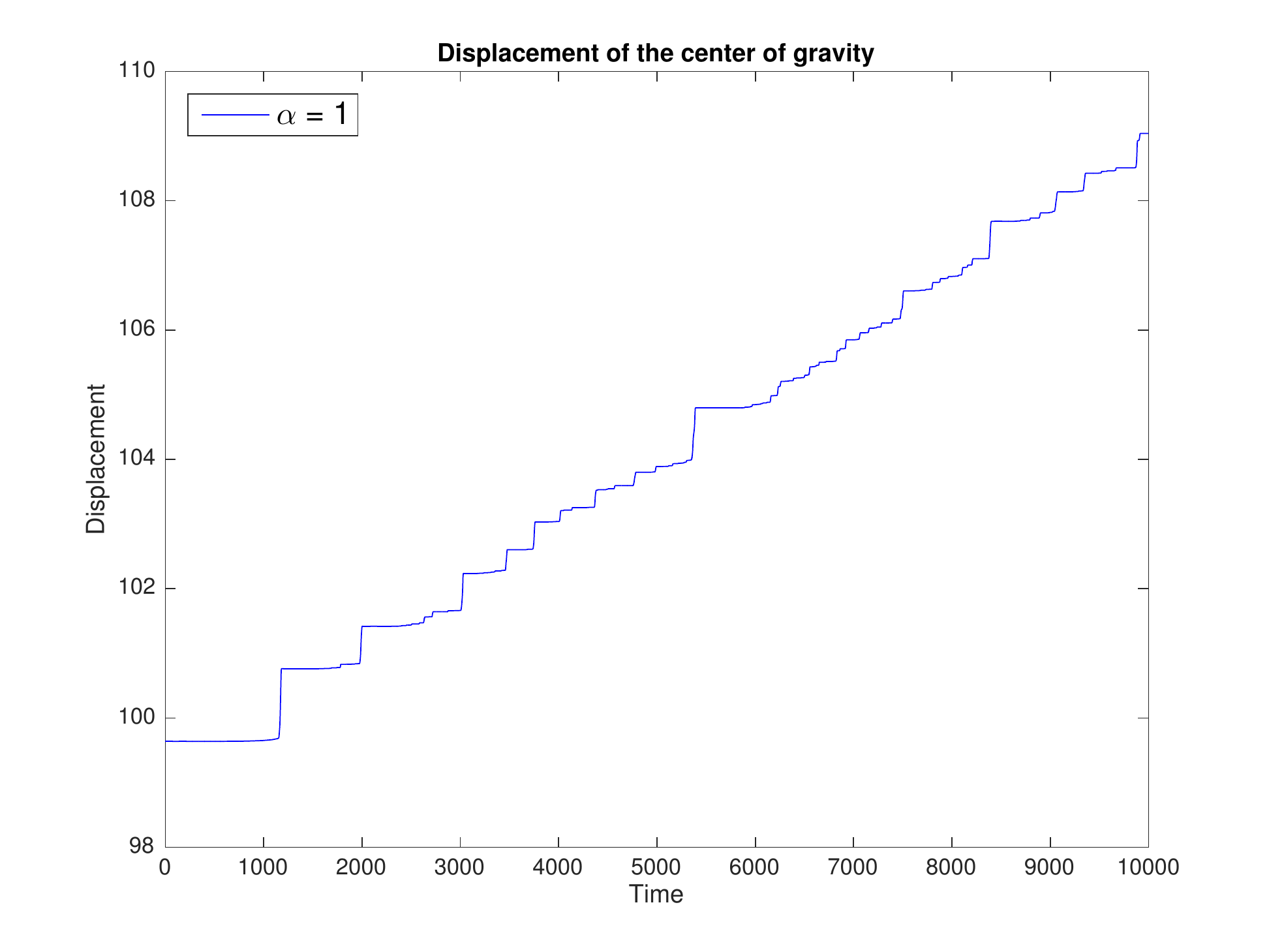}}\quad
\subfloat[][\emph{$\alpha=4$}.]
{\includegraphics[width=.487\columnwidth]{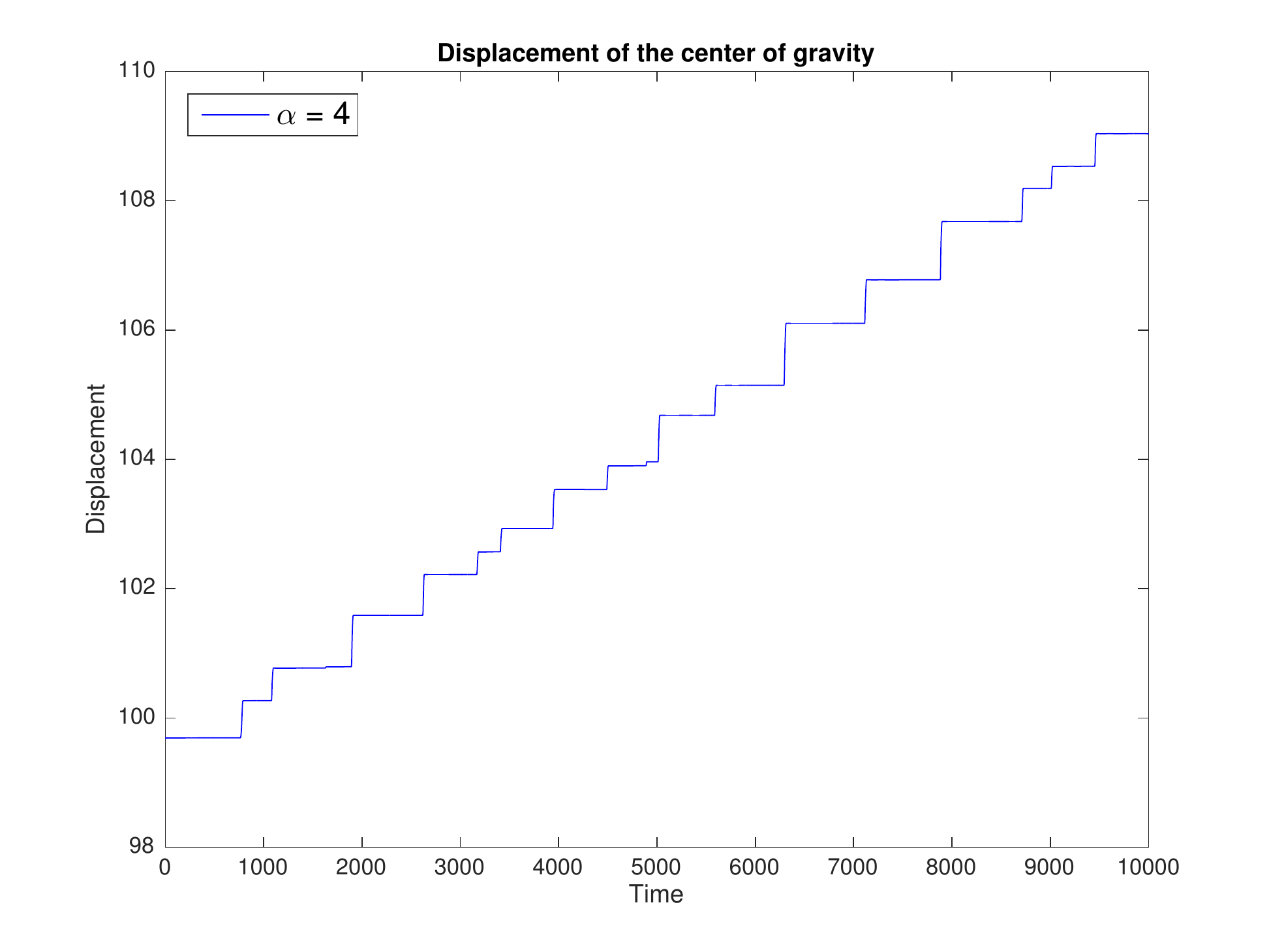}}\quad
\caption{Displacement of the location of the center of gravity: comparison between $\alpha=1$ and $\alpha=4$.}
\label{fig:cdm}
\end{figure}

Now we want to discuss the importance of the parameter $\alpha$ in terms of how much it can affect the results, by pointing out another interesting outcome of the $\alpha$-dependence. We said that the lower the value the higher the dissipation: as a result, in agreement with \cite{CarlLang89}, we notice that in each earthquake the displacements become smaller. By assuming the smallest value employable, namely $\alpha=1$, we can provide a qualitative proof of this property. As in \cite{Saito2007} we consider the displacement of the location of the center of gravity during the time period $[0,10^4]$. Our purpose consists in making a comparison between the cases $\alpha=1$ and $\alpha=4$. Fig.~\ref{fig:cdm} shows the results and indicates that when $\alpha$ is set to unity, the displacements are effectively smaller.

\begin{figure}[!h]
\centering
{\includegraphics[width=.60\columnwidth]{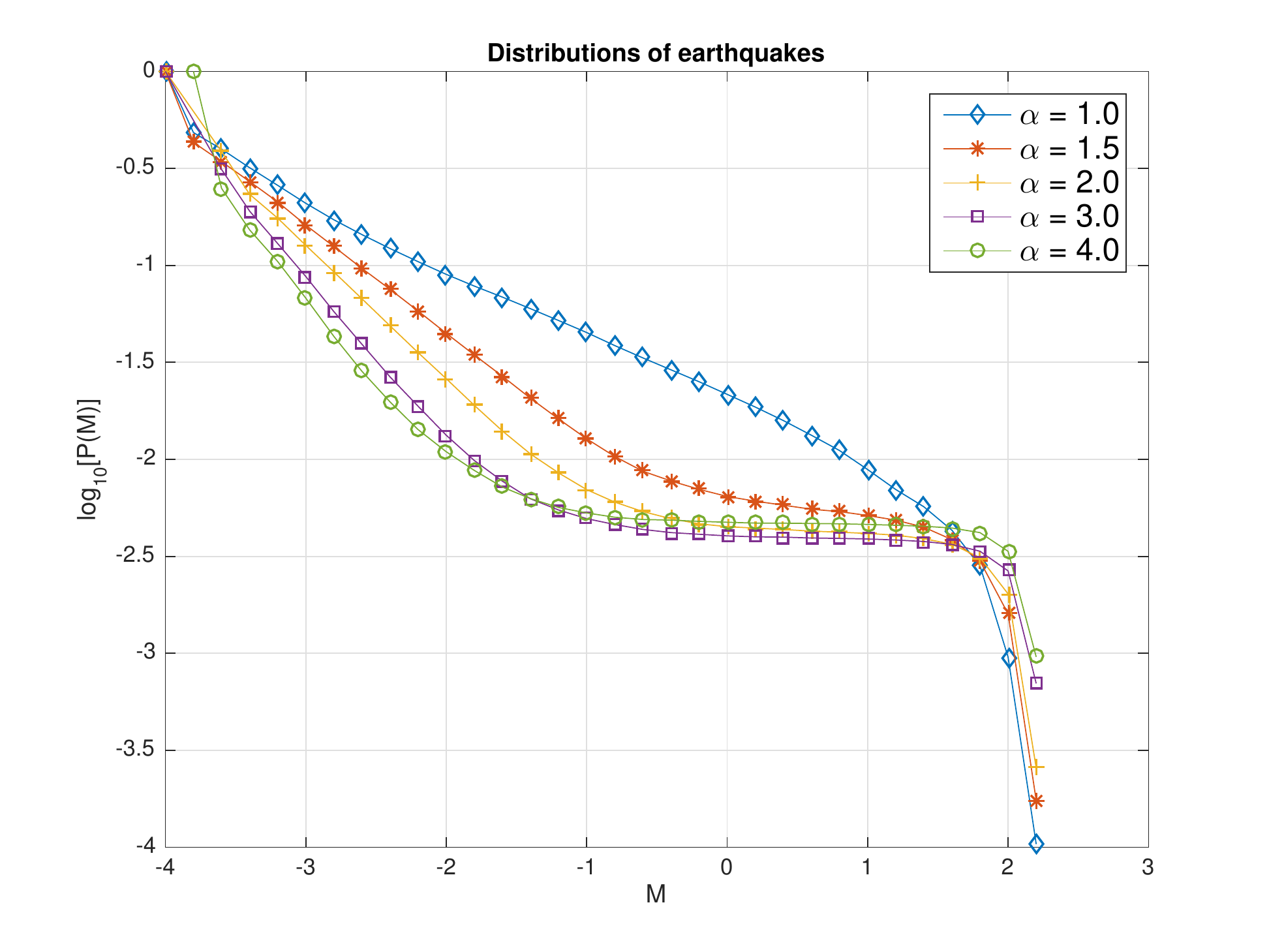}}\quad
\caption{Earthquakes distribution for a system with $200$ blocks: graph of the function $M\mapsto \log_{10}[P(M)]$
with $\alpha \in \{1,1.5,2,3,4\}$. The other parameters used are listed in Table~\ref{tab:param2}.}
\label{fig:TotLog}
\end{figure}

Let us proceed by investigating more accurately how the distribution of earthquakes changes when $\alpha$ increases. In Fig.~\ref{fig:TotLog} some distributions are plotted by assuming $\alpha \in \{1,1.5,2,3,4\}$, so that we have added three different $\alpha$-values beyond those already analysed. These values would be equally spaced if $\alpha=1.5$  was not considered, but we decided to provide a further value within $[1,2]$ in order to control better the evolution of the graph when $\alpha$ is moving in this range. We are allowed to conclude that the peak structure mentioned above persists when $\alpha\not=1$ and that the slope, where there is a linear behaviour qualitatively in agreement with the Gutenberg-Richter law, is steeper as $\alpha$ increases: it can be deduced simply by noticing that when $\alpha=1.5$ a flatter slope affects the distribution in the middle-small range of magnitude, while this slope becomes steeper whether $\alpha$ increases.

Let us conclude this section by making other comparisons of our results. We defined the magnitude $M$ in~\eqref{magni} by introducing the decimal logarithm. Also in \cite{Saito2007} something like this is performed. However, in order to allow further qualitative pairings with some works mentioned in this paper (for example \cite{CarlLang91intrinsic,XionKikuYama15}), we recast the results described above by using the natural logarithm to define the magnitude. The relationship~\eqref{magni} becomes
\begin{equation*}
M_{1}=\ln\biggl(\sum_{i=1}^{N}\Delta x_{i}\biggr),
\end{equation*}
where $M_{1}$ denotes the magnitude. If the quantity $P(M_{1})$ is defined as we have done for $P(M)$, an analogue of Fig.~\ref{fig:TotLog} is obtainable, qualitatively equivalent. We gave it a try and found a graph very similar to the correspondent plots exhibited in \cite{XionKikuYama15}. We also recognized that in our simulations the data, when $\alpha=1$, crosses the other curves at $M_{1}\simeq 4$, in agreement with the results in \cite{XionKikuYama15}.

\begin{figure} [!h]
\centering
\subfloat[][\emph{$\alpha=1$}.]
{\includegraphics[width=.487\columnwidth]{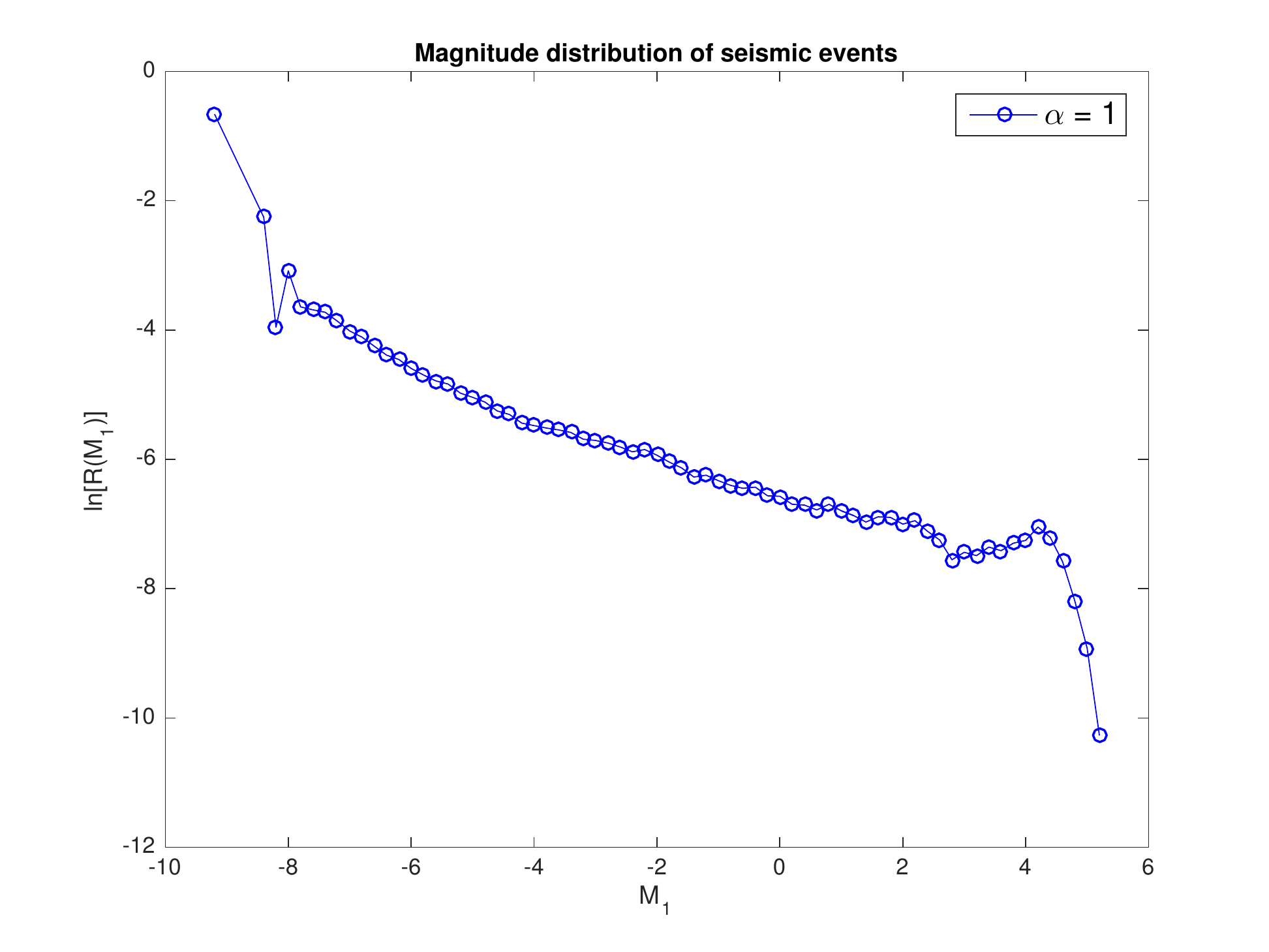}}\quad
\subfloat[][\emph{$\alpha=2$}.]
{\includegraphics[width=.487\columnwidth]{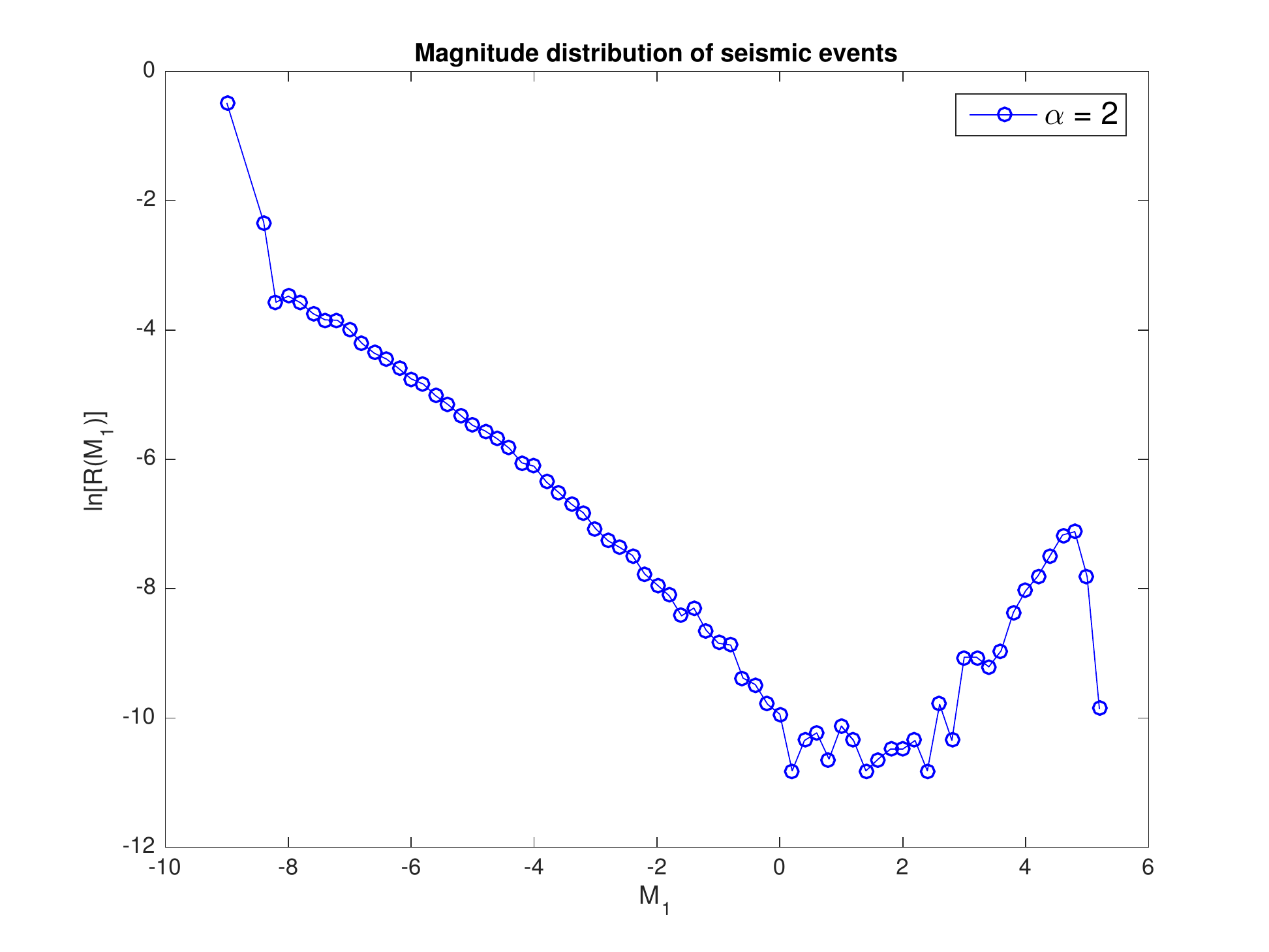}}\quad
\caption{Distributions of magnitude for a system including $200$ blocks: graph of the function $M_{1}\mapsto \ln[R(M_{1})]$;
comparison between $\alpha=1$ and $\alpha=2$.}
\label{fig:CarlsonTipo}
\end{figure}

Finally, if we define another quantity, $R(M_{1})$, as the rate of seismic events with magnitudes equal to $M_{1}$, operatively in a range such as $[M_{1}, M_{1}+dM_{1}]$, we achieve the results shown in Fig.~\ref{fig:CarlsonTipo}.

With the aim of expressing the qualitative behaviour affecting the dynamics when $\alpha\not= 1$, in Fig.~\ref{fig:CarlsonTipo}b we have chosen as example the value $\alpha=2$: it is possible to explain the deviation from the Gutenberg-Richter previously mentioned in terms of a peak structure, by pointing out that in this case large events are too frequent. These conclusions are consistent with those in \cite{CarlLang89,CarlLang91intrinsic,MoriKawa2006}.

\section{Increasing the number of blocks}\label{sec:increasing}
\selectlanguage{english}
As described in Section~\ref{sec:BK}, the Burridge-Knopoff model was originally proposed as a discrete system. As a consequence, despite this kind of representation does not prevent us from achieving interesting and reliable conclusions, we are not allowed to ignore that real earthquakes happen along continuum faults. Therefore, it would be very useful trying to improve our approximation with the aim of providing informations about the continuum limit. In this field we recall some important works \cite{Erick2011,MoriKawa2008}, the first one based on the Dieterich-Ruina friction law, the second one on the velocity-weakening friction proposed by Carlson and Langer \cite{CarlLang89properties,CarlLang89}.

Although we will pursue the issue of deducing a suitable continuum version of the model in a forthcoming work, we propose a preliminary
approach to this goal by considering again the discrete version and by increasing the number of blocks within the system.

Our strategy consists, firstly, in deducing the continuum version starting from the relationship~\eqref{edoi}; secondly, we will understand how to rescale the parameters in order to describe the same configuration by increasing the number of blocks; finally, by analysing both the distance estimates of some magnitude distributions and their qualitative behaviours, we aim at establishing whether the convergence is satisfying or not.
We recall the equation~\eqref{edoi}, that is
\begin{equation}
\label{oneblocktime}
m\ddot{x}_{i}=k_{c}(x_{i+1}-2x_{i}+x_{i-1})+k_{p}(Vt-x_{i})-F(\dot{x}_{i}).
\end{equation} 
It is important remembering that we supposed the blocks to rest on a surface considered as a one-dimensional object. Moreover, the blocks were assumed to be equally spaced with distance $a$. This parameter plays the role of the mesh size within the one-dimensional grid (we indicate with $L$ the grid length) defined by all the blocks and used to provide a discrete representation of one of the sides of a fault. That is why, to obtain a continuum model, whose solution is spatial and time dependent, we have to take the limit $a \rightarrow 0$. Because $a$ does not explicitly appear in~\eqref{oneblocktime}, we firstly divide both sides by this parameter and get   
\begin{equation}
 \label{oneblocktimeb}
\frac{m\ddot{x}_{i}}{a}=ak_{c} \left( \frac{x_{i+1}-2x_{i}+x_{i-1}}{a^2} \right) + \frac{k_{p}(Vt-x_{i})}{a} - \frac{F(\dot{x}_{i})}{a}.
\end{equation}
It is useful to notice that the quantity $L$ can be related to the number of blocks, $N$, by
\begin{equation}
\label{LatticeBlocks}
L=Na.
\end{equation}
This means that we can equivalently assume $a \rightarrow 0$ or rather $N \rightarrow \infty$ to take the continuum limit. This double way of thinking is very advantageous: on the one hand, it allows us to approximate operatively the continuum limit by increasing the number of blocks; on the other, it makes easy obtaining, as limit of~\eqref{oneblocktimeb}, a forced wave equation with damping and friction for some macroscopic displacement 
$u=u(x,t)$, namely
\begin{equation}
\label{continuum}
M\partial_{tt} u = K_{c}\partial_{xx} u+K_{p}(Vt-u)-\phi(\partial_{t} u). 
\end{equation}
Due to the fact that $L$ is a constant, we have simply imposed $L=1$ in~\eqref{LatticeBlocks}.
As a consequence, the quantities involved in~\eqref{continuum} are definable as in Table~\ref{tab:param3}. 
\begin{table}[!h]
\caption{Quantities involved in~\eqref{continuum}. }
\label{tab:param3}
\centering
\begin{tabular}{*{3}{clr}}
\toprule
\multicolumn{3}{c}%
{\textbf{Matching}}   \\   
\midrule
$Nm^{(N)}$ & $=$ & $M$ \\
$N^{-1}k^{(N)}_{c}$ & $=$ & $K_{c}$ \\
$Nk^{(N)}_{p}$ & $=$ & $ K_{p} $\\ 
$NF^{(N)}(\dot{x}_{i})$ & $=$ & $\phi(\partial_{t} u)$\\ 
\bottomrule
\end{tabular}
\end{table}
We have included the superscript $N$ to identify the belonging to the discrete system involving $N$ blocks. 

Because in this paper we do not simulate the Burridge-Knopoff model in a continuum form as in~\eqref{continuum}, we do not provide a more explicit definition for the friction law $ \phi(\partial_{t} u)$ and the other parameters. Our goal consists solely in understanding how to rescale the quantities appearing in~\eqref{oneblocktime} taking advantage of~\eqref{continuum}, in order to investigate the continuum limit by using the discrete version. For instance, let us consider the parameter $M$: it can be easily associated to the total mass belonging to a part of the fault. The best approximation of this total mass is performed in the continuum version but it can be approximated also with an ideal chain including a number of blocks tending to infinity, whose mass tends to be infinitesimal, within the discrete model. This point of view allows us to approach a hypothetical continuum limit by increasing the number of blocks and updating the parameters $m^{(N)}, k^{(N)}_{c},k^{(N)}_{p}$ and the friction law appropriately. More specifically, starting with $N$ blocks, if we choose to double the blocks, we would double $k^{(N)}_{c}$ and halve the remaining  $N$-dependent quantities appearing in Table~\ref{tab:param3}.

\begin{figure}[!h]
\centering
{\includegraphics[width=.60\columnwidth]{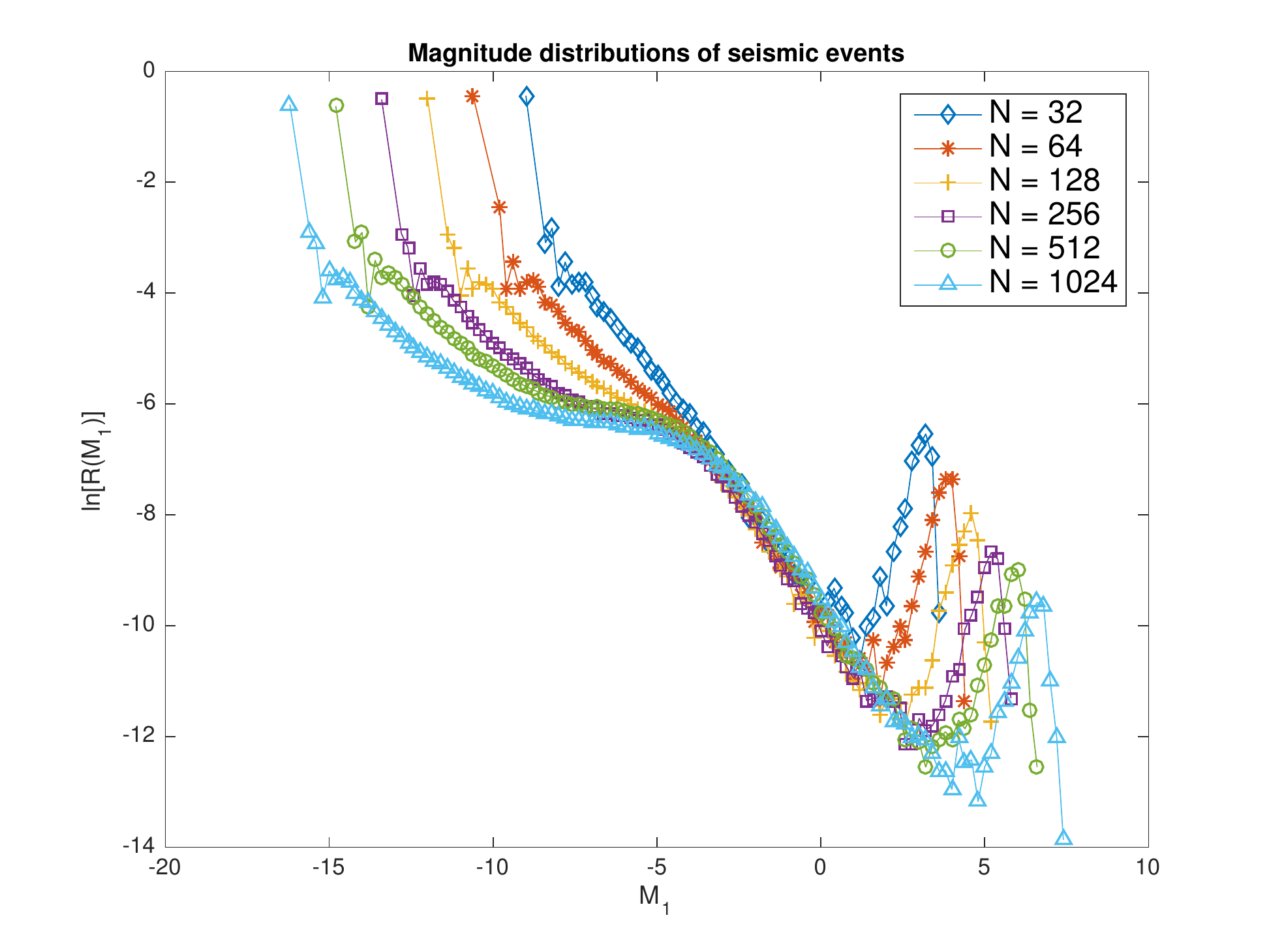}}\quad
\caption{Magnitude distributions of earthquakes to investigate the continuum limit.
The number of blocks starts from $N=32$ and ends at $N=1024$.}
\label{fig:Continuum}
\end{figure}

Let us start with $N=32$, firstly adopting the parameters listed in Table~\ref{tab:param4}. We proceed by doubling five times the number of blocks and updating the  $N$-dependent quantities as explained above.
As regards the friction law, we simply update it by changing the $F^{(N)}_{0}$ value. 
\begin{table}[!h]
\caption{Values for the quantities involved in~\eqref{continuum} when $N=32$.}
\label{tab:param4}
\centering
\begin{tabular}{*{7}{c}}
\toprule
\multicolumn{7}{c}%
{\textbf{Parameters}}  \\   
\midrule
$\mathbf{m^{(N)}}$ &$\mathbf{k^{(N)}_{p}}$ &$\mathbf{k^{(N)}_{c}}$ &$\mathbf{V}$ &$\mathbf{F^{(N)}_{0}}$  &$\mathbf{\sigma}$  & $\mathbf{\alpha}$\\  
\midrule  
  $1$      & $1$       & $100$      & $0.001$ & $1$ & $0.01$ & $2$  \\
\bottomrule
\end{tabular}
\end{table}
We choose to represent the magnitude distributions by using the quantity $R(M_{1})$ defined in Section~\ref{sec:algorithm}, as performed in \cite{MoriKawa2008}. Because of the five doubling of $N$, simulations with $N \in \{32,64,128,256,512,1024\}$ are performed. We also introduce an additional parameter, $n$, to quickly recall the number of doubling related to the $N$ blocks involved: for instance, $N=32$ corresponds to $n=0$, $N=64$ to $n=1$ and so on, until $n=5$. The results are shown in Fig.~\ref{fig:Continuum}.
Finally, in order to provide some distance estimates among the data, we compute both the norms $\Vert \bullet \Vert_{2} $ and $\Vert \bullet \Vert_{\infty}$, by using the magnitude distributions, $f_{n}$, restricted to the set of common magnitudes. The biggest set corresponds to the magnitude distribution provided by $N=32$, which is the biggest intersection: we select almost the whole set, except for the first value that, as we pointed out in Section~\ref{sec:algorithm}, can be partially influenced by the discreteness. The results are reported in Table~\ref{tab:estimates}.
\begin{table}[!h]
\caption{Distance estimates.}
\label{tab:estimates}
\centering
\begin{tabular}{*{3}{c}}
\toprule
 & $\Vert \bullet \Vert_{2} $  & $\Vert \bullet \Vert_{\infty}$ \\
\midrule
 $f_{0}-f_{1}$  & $7.6788$  & $2.6056$ \\
 $f_{1}-f_{2}$  & $6.1455$  & $2.5403$ \\
 $f_{2}-f_{3}$  & $3.4142$  & $1.8730$ \\
 $f_{3}-f_{4}$  & $2.1490$  & $0.6726$ \\
 $f_{4}-f_{5}$  & $1.7111$  & $0.5737$ \\
\bottomrule
\end{tabular}
\end{table}

Although the distance estimates seem to be very comforting, because of the constant decreasing provided by Table~\ref{tab:estimates}, we are not allowed to recognize a convergence to an asymptotic law, in agreement with \cite{MoriKawa2008}. Of course, the distributions tend to be narrower as the number of blocks increases but at the same time these laws change structure, more specifically becoming wider, as can be deduced by analysing the plots in Fig.~\ref{fig:Continuum}. This phenomenon can be observed mostly close to the edges of the graphs: on the left side magnitude becomes smaller, on the right one it gets larger, as $n$ increases. All these results are consistent with the work by Mori and Kawamura \cite{MoriKawa2008}. As concerns the smallest magnitudes this same work suggests that the tendency of being smaller is explainable in terms of infinitesimal earthquakes potentially occurring in the continuum limit.   

\section{Discussion and perspectives}\label{sec:disc}
\selectlanguage{english}
The Predictor-Corrector numerical strategy has been performed in this paper for the simulation of the Burridge-Knopoff model
and it has been proven to be effective for reproducing real behaviors.
This method allows us to employ a total explicit procedure, without solving a nonlinear system at every step as in the case of an implicit method, but keeping a good quality in terms of simulations. Our results, indeed, are consistent with previous works on the same model: by making a comparison with the Gutenberg-Richter law, it is possible to recognize a deviation at large magnitudes for $\alpha>1$; for $\alpha=1$ a qualitative behaviour in agreement with the empirical law is noticeable, although with a flatter slope. We have also studied the dynamics for the simplest cases, by involving few blocks, in order to discuss the main features related to the stick-slip phenomenon. Finally, we have investigated the continuum limit by rescaling the quantities included within the discrete model and increasing the number of blocks. Our results, in agreement with \cite{MoriKawa2008}, point out that although the magnitude distributions tend to be narrower, a full convergence to an asymptotic form is not achieved.

Future works is aimed at performing tests by using an optimized implicit method which allows to consider a very huge number of blocks, thus bypassing the problem of the size of the system. In this sense the attempt proposed in \cite{XionKikuYama15} paved the way. The perspective of an implicit method is very useful also because it does not prevent from using larger step size which would be inadvisable within an explicit strategy.
Moreover, it would be possible to improve the computational effectiveness by adopting the parallel computation paradigm.

Of course, as regards the model and the physical contributions involved, settings aimed at reproducing more realistic trends must be employed, maybe by acting on the choice of the friction law or the elastic forces.  

\section*{Acknowledgements}
We sincerely thank Dr. Xiaogang Xiong for the fruitful exchange of informations on various topics related to the Burridge-Knopoff model.
We also like to thank Prof. Nicola Guglielmi for his suggestions about the possible developments related to the numerical issues,
and Dr. Chiara Simeoni for her contribution especially in the writing process of the manuscript.

\end{document}